\documentclass[11pt]{amsart}
\usepackage[english]{babel}
\usepackage[latin1]{inputenc}
\usepackage{amsmath}
\usepackage{amsthm}
\usepackage{amsfonts}
\usepackage{amssymb, latexsym, graphics, showidx}
\usepackage{graphicx}

\numberwithin{equation}{section}

\newtheorem{thm}{Theorem}[section]
\newtheorem{rem}[thm]{Remark}

\newtheorem{prop}[thm]{Proposition}
\newtheorem{lem}[thm]{Lemma}
\newtheorem*{ack}{Acknowledgments}

\renewcommand{\dim}{\noindent\textbf{Proof.} }
\newcommand{\dims}{\noindent\textbf{Proof} }
\newcommand{\finedim}{{\unskip\nobreak\hfil\penalty50
   \hskip2em\hbox{}\nobreak\hfil\mbox{$\Box$ \qquad}
   \parfillskip=0pt \finalhyphendemerits=0\par\medskip}}
\newcommand{\R}{\mathbb{R}}
\newcommand{\N}{\mathbb{N}}
\newcommand{\Z}{\mathbb{Z}}

\newcommand{\lam}{\lambda}
\newcommand{\al}{\alpha}
\newcommand{\lams}{\overline{\lambda}}

\newcommand{\xs}{\overline{x}}
\newcommand{\ys}{\overline{y}}

\newcommand{\ts}{\overline{t}}
\newcommand{\s}{\overline{s}}
\newcommand{\Xs}{\overline{X}}
\newcommand{\Ys}{\overline{Y}}

\newcommand{\ep}{\epsilon}

\newcommand{\beq}{\begin{equation}}
\newcommand{\eeq}{\end{equation}}
\newcommand{\beqs}{\begin{equation*}}
\newcommand{\eeqs}{\end{equation*}}
\newcommand{\beqa}{\begin{eqnarray}}
\newcommand{\eeqa}{\end{eqnarray}}
\newcommand{\beqas}{\begin{eqnarray*}}
\newcommand{\eeqas}{\end{eqnarray*}}

\begin{document}

\title[Homogenization of first order equations]{Homogenization of first order equations \\with $u/\epsilon$-periodic Hamiltonian:\\ Rate of convergence as $\epsilon\to 0$\\ and numerical approximation of the effective Hamiltonian. }

\author{Yves Achdou}
\address{ UFR Math{\'e}matiques, Universit{\'e} Paris Diderot, Case 7012,
75251 Paris Cedex 05, France \\
and Laboratoire Jacques-Louis Lions, Universit{\'e} Paris 6, 75252 Paris Cedex 05}
\email{
achdou@math.jussieu.fr
}

\author{Stefania Patrizi}
\address{SAPIENZA Universit{\`a} di Roma, Dipartimento di Matematica,
Piazzale A.~Moro 2, I-00185 Roma, Italy}
\email{patrizi@mat.uniroma1.it}

\maketitle

\begin{abstract}
  We  consider homogenization problems for first order Hamilton-Jacobi equations with $u^\epsilon/\epsilon$ periodic dependence, recently introduced by C. Imbert and R. Monneau, and also studied by G. Barles: this unusual dependence leads to a nonstandard cell problems.
We study the rate of convergence of the solution to the solution of the homogenized problem when the parameter $\epsilon$ tends to $0$. We obtain the same rates as those obtained by I. Capuzzo Dolcetta and H. Ishii for the more usual homogenization problems without the dependence in $u^\epsilon/\epsilon$. In a second part, we study Eulerian schemes for the approximation of the cell problems.
We prove that  when the grid steps tend to zero, the approximation of the effective Hamiltonian converges to the effective Hamiltonian.
\end{abstract}

\section{Introduction}
We consider homogenization problems for first order Hamilton-Jacobi equations with $u^\epsilon/\epsilon$ periodic dependence, namely
\beq\label{uep} \left\{%
\begin{array}{ll}
    u^\ep_t+H\left(\frac{t}{\ep},\frac{x}{\ep},\frac{u^\ep}{\ep},Du^\ep\right)=0, & (t,x)\in(0,+\infty)\times\R^N, \\
    u^\ep(0,x)=u_0(x), & x\in\R^N\\
\end{array}%
\right. \eeq
with the following assumptions on the Hamiltonian $H$:
\renewcommand{\labelenumi}{(H\arabic{enumi})}
\begin{enumerate}
\item Periodicity: for any $(t,x,u,p)\in \R\times\R^N\times\R\times\R^N$
$$ H(t+1,x+k,u+1,p)=H(t,x,u,p)\quad \text{for any
}k\in\Z^N;$$
 \item Regularity: $H:\R\times\R^N\times\R\times\R^N\rightarrow\R$ is Lipschitz
 continuous and there exists a constant $C_1>0$ such that, for
 almost every $(t,x,u,p)\in \R\times\R^N\times\R\times\R^N$
 $$|D_{(t,x)} H(t,x,u,p)|\leq C_1(1+|p|),\quad |D_u H(t,x,u,p)|\leq
 C_1,\quad|D_p H(t,x,u,p)|\leq C_1;$$

 \item $H(t,x,u,p)\rightarrow+\infty$ as $|p|\rightarrow+\infty$ uniformly for $(t,x,u)\in
 \R\times\R^N\times\R$;

 \item There exists a constant $C$ such that for almost every $(t,x,u,p)\in \R\times\R^N\times\R\times\R^N$
 $$|D_p H(t,x,u,p)\cdot p-H(t,x,u,p)|\leq C.$$
 \end{enumerate}

Problem \eqref{uep} with $H$ independent of $t$  was introduced by Imbert and Monneau \cite{im} as a simplified  model for dislocation dynamics in material science.
The complete model is introduced in \cite{imr} and leads to nonlocal first order  equations of the type
\begin{displaymath}
    u^\ep_t +\left(c(\frac x \ep)+ M^\epsilon(\frac{u^\ep}{\ep})\right) | Du^\ep|  +    H\left(\frac{u^\ep}{\ep},Du^\ep\right)=0
\end{displaymath}
where $M^\epsilon$ is a nonlocal jump operator and $c$ is a periodic velocity. In the latter model, the level sets of the solution $u^\ep$ describe dislocations.

Going back to \eqref{uep}, it was proved in \cite{im} that, with $H$ independent of $t$,
\begin{itemize}
\item under  assumptions (H1) and (H2), there exists a unique bounded continuous viscosity solution of \eqref{uep};
\item under assumptions (H1)-(H3),  the limit $u^0$ of $u^\ep$ as $\ep\to0$ exists and it is the unique bounded continuous solution of the homogenized problem
\beq\label{u0} \left\{%
\begin{array}{ll}
    u^0_t+\overline{H}(Du^0)=0, & (t,x)\in(0,+\infty)\times\R^N, \\
    u^0(0,x)=u_0(x), & x\in\R^N,\\
\end{array}%
\right. \eeq
where the effective Hamiltonian $\overline{H}$ is uniquely defined by the long time behavior of the solution of
\beq\label{v} \left\{%
\begin{array}{ll}
\lam=v_t+H(x,-\lam t+p\cdot x+v,p+Dv), & (t,x)\in(0,+\infty)\times\R^N, \\
    v(0,x)=0, & x\in\R^N.\\
\end{array}%
\right. \eeq
\end{itemize}
More precisely, we have the following theorem
\begin{thm}[Imbert-Monneau, \cite{im}]\label{sec:th1}
Let $H$ be independent of $t$. Assume  (H1)-(H3) and $u_0\in
W^{1,\infty}(\R^N)$. Then, as $\ep\rightarrow0$, the sequence
$u^\ep$ converges locally uniformly in $(0,+\infty)\times\R^N$ to
the solution $u^0$ of \eqref{u0}, where, for any $p\in\R^N$
$\overline{H}(p)$ is defined as the unique number $\lam$ for which
there exists a bounded continuous viscosity solution of \eqref{v}.
Moreover $\overline{H}:\R^N\rightarrow\R$ is continuous and
satisfies the coercivity property
$$\overline{H}(p)\rightarrow +\infty\quad\text{as
}|p|\rightarrow+\infty.$$ \end{thm}
The proof in \cite{im} is rather involved: it uses a {\it twisted} perturbed test function for a higher dimensional problem posed in $\R\times \R^{N}\times \R$.\\
Under the additional assumption (H4), an easier proof of Theorem  \ref{sec:th1} was given by Barles, \cite{b},
 as a byproduct of a general result on the homogenization of Hamilton-Jacobi equations with non-coercive Hamiltonians.
\begin{rem}\label{rem1}{\em The hypothesis (H4) which was not used in \cite{im} guarantees the existence of a function
$H_\infty$ such that
$$H_\infty(t,x,u,p)=\lim_{s\rightarrow0^+}sH(t,x,u,s^{-1}p).$$
Moreover $H_\infty$ satisfies (H1)-(H3).}
\end{rem}
In \cite{b}, thanks to assumption (H4), the equation for $u^\epsilon$ is interpreted as an equation for the motion of a graph: indeed, following \cite{b}, for $t\in\R$, $(x,y)\in \R^{N+1}$,
$(p_x,p_y)\in \R^{N+1}$, let us introduce the non-coercive
Hamiltonian $F$ defined by
\beq\label{F} F(t,x,y,p_x,p_y)=\left\{%
\begin{array}{ll}
    |p_y|H(t,x,y,|p_y|^{-1}p_x), & \hbox{if }p_y\neq 0, \\
    H_\infty(t,x,y,p_x), & \hbox{otherwise.} \\
\end{array}%
\right. \eeq
The function $U^\ep(t,x,y):=u^\ep(t,x)-y$ satisfies
 \beq\label{Uep} \left\{%
\begin{array}{ll}
    U^\ep_t+F\left(\frac{t}{\ep},\frac{x}{\ep},\frac{U^\ep+y}{\ep},D_xU^\ep,D_yU^\ep\right)=0, & (t,x,y)\in(0,+\infty)\times\R^{N+1}, \\
    U^\ep(0,x,y)=u_0(x)-y, & (x,y)\in\R^{N+1}.\\
\end{array}%
\right. \eeq
 In \cite{b} Barles proves that the sequence
$U^\ep$ converges to the solution $U^0$ of the following problem
\beq\label{U0}\left\{%
\begin{array}{ll}
    U^0_t+\overline{F}(D_x U^0,D_yU^0)=0, & (t,x,y)\in(0,+\infty)\times\R^{N+1}, \\
    U^0(0,x,y)=u_0(x)-y, & (x,y)\in\R^{N+1},\\
\end{array}%
\right. \eeq where for  $(p_x,p_y)\in\R^{N+1}$,
$\overline{F}(p_x,p_y)$ is the unique number $\lam$ for which the
cell problem \beq\label{V}
V_t+F(t,x,y,p_x+D_xV,p_y+D_yV)=\lam\quad\text{in
}\R\times\R^{N+1}.\eeq
 admits bounded sub and supersolutions. This result makes it possible to solve the homogenization
problem for \eqref{uep}:
\begin{thm}[Barles, \cite{b}]\label{convthmbarles} Assume (H1)-(H4). Then the sequence
$u^\ep$ converges locally uniformly in $(0,+\infty)\times\R^{N}$
to the solution $u^0$ of \eqref{u0}. The function
$\overline{H}(p)$ in \eqref{u0} can be cha\-racte\-ri\-zed as
follows: $\overline{H}(p)=\overline{F}(p,-1)$, where, for any
$(p_x,p_y)\in\R^{N+1}$, $\overline{F}(p_x,p_y)$ is the unique
number $\lam$ for which the equation \eqref{V} admits bounded sub
and supersolutions in $\R\times\R^{N+1}$.
\end{thm}
An important step in the proof of Theorem \ref{convthmbarles} consists of homogenizing the non-coercive level-set equation satisfied by $1\!\!1 _{\{U^\epsilon\ge 0\}}$.
\\
In this paper, we  tackle two questions:
\begin{itemize}
\item Is it possible to estimate the rate of convergence of $u^\ep$ to $u^0$ when $\epsilon\to 0$?
\item Is is possible to approximate numerically the effective Hamiltonian?
\end{itemize}
The first question was answered by Capuzzo Dolcetta and Ishii, \cite{ci} for a more classical homogenization problem: the estimate $\|u^\epsilon-u^0\|_{\infty} \le C\epsilon^{\frac 1 3}$ was obtained for Hamilton-Jacobi equations of the type
\begin{displaymath}
  u^\epsilon +H\left(x,\frac x \epsilon, u^\epsilon\right)=0,
\end{displaymath}
where $(x,y,p)\to H(x,y,p)$ is a coercive Hamiltonian,  uniformly Lipschitz continuous  for $|p|$ bounded and  periodic with respect to $y$; moreover, if $ H(x,y,p)$ does not depend on $x$, then the convergence is linear in $\epsilon$. We will show that in the present case, it is possible to obtain the same rates of convergence as $\epsilon\to 0$ by adapting the proof in \cite{ci} using the arguments contained in \cite{b}. Our main result on this topic is Theorem~\ref{main} in \S~\ref{sec:an-estimate-rate}. The main idea is to approximate $U^\epsilon$ (with an error smaller than $\epsilon$) by a discontinuous function $\widetilde U^\epsilon$ which takes integer values where $U^\ep$ has noninteger values and  which is a discontinuous viscosity solution of
\[  \widetilde U^\ep_t+F\left(\frac{t}{\ep},\frac{x}{\ep},\frac{y}{\ep},D_x \widetilde U^\ep,D_y\widetilde U^\ep\right)=0, \quad\quad (t,x,y)\in(0,+\infty)\times\R^{N+1}.\]
The latter equation  has to be compared with \eqref{Uep}.
This approximation $ \widetilde U^\ep$ is obtained as the limit as $\delta\to 0$ of $\phi_\delta(U^\ep)$ where $(\phi_\delta)_{\delta}$ is a sequence of increasing functions. The method of Capuzzo Dolcetta and  Ishii \cite{ci} can then be applied to $\widetilde U^\epsilon$.
\\
The second question was studied in \cite{acc} for equation
\begin{displaymath}
  u^\epsilon +H\left(\frac x \epsilon, u^\epsilon\right)=0,
\end{displaymath}
where $(y,p)\to H(y,p)$ is a coercive Hamiltonian,  uniformly
Lipschitz continuous  for $|p|$ bounded and  periodic with respect
to $y$; in this article, a complete numerical method for solving
the homogenized problem was studied, including as a main step the
approximation of the effective Hamiltonian by solving discrete
cell problems. Error estimates were proved. Here, we will  study
the approximation of the cell problem (\ref{V}) by Eulerian
schemes in the discrete torus. We have prefered to study the
approximation of the noncoercive $N+2$ dimensional problem
(\ref{V})  rather than that of the coercive $N+1$ dimensional
problem (\ref{v})   because the solution of (\ref{v}) may not be
periodic. In \S~\ref{sec:appr-effect-hamilt}, we prove
Theorem~\ref{ergodicapprox}, the discrete analogue of the
ergodicity Theorems in \cite{b}, i.e. that there exists a unique
real number $\lambda_h^{\Delta t}$ such that the discrete analogue
of \eqref{V} has a solution. The arguments in the proof are the
discrete counterparts of those in \cite{b}. Then, we prove
Proposition~\ref{convappham}, which states that the discrete
effective Hamiltonian converges to the effective Hamiltonian when
the grid step of the discrete cell problem tends to zero.
\\
To summarize, the paper is organized as follows:
Section \ref{sec:an-estimate-rate} is devoted to finding estimates on the rate of convergence as $\epsilon\to 0$.
Section \ref{sec:appr-effect-hamilt} is devoted to the numerical approximation of the effective Hamiltonian by Eulerian schemes.
Finally, we present some numerical tests in  Section \ref{sec:numerical-tests}.

\section{An estimate on  the rate of convergence when  $\epsilon\rightarrow
0$} \label{sec:an-estimate-rate}
This section is devoted to the estimate of the rate of the
uniform convergence of the solutions of \eqref{uep} to the
solution of the equation \eqref{u0} in term of $\ep$.
\subsection{The main result}
\begin{thm}\label{main} Assume (H1)-(H4) and $u_0\in W^{1,\infty}(\R^N)$. Let
$u^\ep$ and $u^0$ be respectively the viscosity solutions of
\eqref{uep} and \eqref{u0}. Then there exists a constant $C$,
independent of $\ep\in(0,1)$, such that for any $T>0$
\beq\label{main1} \sup_{[0,T]\times\R^N}|u^\ep(t,x)-u^0(t,x)|\leq
C e^T\ep^\frac{1}{3}.\eeq  If $u_0$ is affine then
\beq\label{main2}\sup_{\R^+\times\R^N}|u^\ep(t,x)-u^0(t,x)|\leq C
\ep.\eeq
\end{thm}

\subsection{Preliminary results} In this section we recall some
results that will be used later to obtain error estimates.

 The assumptions (H1)-(H4) on $H$
guarantee that $F$ satisfies
\renewcommand{\labelenumi}{(F\arabic{enumi})}
\begin{enumerate}
\item Periodicity: for any $(t,x,y,p_x,p_y)\in \R\times\R^{N+1}\times\R^{N+1}$
$$ F(t+1,x+k,y+1,p_x,p_y)=F(t,x,y,p_x,p_y)\quad \text{for any
}k\in\Z^N;$$
 \item Regularity: $F:\R\times\R^{N+1}\times\R^{N+1}\rightarrow\R$ is Lipschitz
 continuous and there exists a constant $C_1>0$ such that, for
 almost every $(t,x,y,p_x,p_y)\in \R\times\R^{N+1}\times\R^{N+1}$
 $$|D_{(t,x)} F(t,x,y,p_x,p_y)|\leq C_1(|p_x|+|p_y|),\; |D_y F(t,x,y,p_x,p_y)|\leq
 C_1|p_y|,$$ $$|D_{(p_x,p_y)} F(t,x,y,p_x,p_y)|\leq C_1;$$

 \item Coercivity: $F(t,x,y,p_x,p_y)\rightarrow+\infty$ as $|p_x|\rightarrow+\infty$ uniformly for $(t,x,y)\in
 \R\times\R^{N+1}$, $|p_y|\leq R$, for any $R>0$;
\end{enumerate}
Remark that $F(t,x,y,0,0)=0$. This and (F2) imply that for every
$(t,x,y,p_x,p_y)\in \R\times\R^{N+1}\times\R^{N+1}$
 \beq\label{F4} |F(t,x,y,p_x,p_y)|\leq C_1(|p_x|+|p_y|).
 \eeq

 Moreover, by construction, $F$ satisfies the "geometrical"
 assumption
 \begin{enumerate}
\item[(F4)] For any $(t,x,y,p_x,p_y)\in \R\times\R^{N+1}\times\R^{N+1}$ and
any $\lam>0$,
$$F(t,x,y,\lam p_x,\lam p_y)=\lam F(t,x,y,p_x,p_y).$$
\end{enumerate}
Assumption (F4) guarantees that \eqref{Uep} is invariant by
any nondecreasing change $U\rightarrow\varphi(U)$, see \cite{cgg} and
\cite{es}, i.e., any function $V=\varphi(U^\ep)$, with $\varphi$
nondecreasing is solution of
 \beqs \left\{%
\begin{array}{ll}
    V_t+F\left(\frac{t}{\ep},\frac{x}{\ep},\frac{U^\ep+y}{\ep},D_xV,D_yV\right)=0, & (t,x,y)\in(0,+\infty)\times\R^{N+1}, \\
    V(0,x,y)=\varphi(u_0(x)-y), & (x,y)\in\R^{N+1}.\\
\end{array}%
\right. \eeqs Finally, note that  (F3) and (F4) imply the
existence of a positive constant $C_2$ such that
\beq\label{Fcoercive}F(t,x,y,p_x,0)\geq C_2|p_x| \quad\text{for
all }(t,x,y,p_x)\in\R\times\R^{N+1}\times\R^N.\eeq

In \cite{b},
in order to construct sub and supersolutions of \eqref{V},
Barles introduces for $\al>0$ the auxiliary equation
\beq\label{Wal} W^\al_t+F(t,x,y,p_x+D_xW^\al,p_y+D_yW^\al)+\al
W^\al=0,\quad (t,x,y)\in \R\times\R^{N+1},\eeq
with $F$ defined by \eqref{F}, and shows that
% and shows that under
% suitable assumptions on $F$,
 if (H1)-(H4)
hold true, then \eqref{Wal} admits a unique continuous periodic
viscosity solution. Moreover the limit of $\al W^\al(t,x,y)$ as
$\al\rightarrow0^+$ does not depend on $(t,x,y)$ and the
half-relaxed limits of $W^\al-\min W^\al$ provide a bounded
subsolution and a bounded supersolution of \eqref{V}, with
$\lam=-\lim_{\al\rightarrow0^+}\al W^\al(t,x,y)$. We use the notation
$P=(p_x,p_y)\in\R^{N+1}$ and  $W^\al(x,y,P)$ for the unique solution
of \eqref{Wal}. We have the following proposition:
\begin{prop}[Barles, \cite{b}]\label{barlesprop}
For any $(t,x,y,P)\in \R\times\R^{N+1}\times\R^{N+1}$,
$P=(p_x,p_y)$, the following estimates hold
\begin{itemize}
    \item [(i)] $\min_{(t,x,y)\in\R\times\R^{N+1}} -F(t,x,y,P)\leq \al
    W^\al(t,x,y,P)\leq\max_{(t,x,y)\in\R\times\R^{N+1}}
    -F(t,x,y,P)$;
    \item [(ii)] There exists a constant $K_1>0$ depending on $\|F(t,x,y,p_x,p_y)\|_\infty$ and $C_2$
    such that
    $$ \max_{\R\times\R^{N+1}} W^\al-\min_{\R\times\R^{N+1}} W^\al\leq K_1.$$
    \end{itemize}
\end{prop}
Further properties of  $W^\al(x,y,P)$ are given in the following
lemma:

\begin{lem}\label{wproplem}For any $(t,x,y,P)\in \R\times\R^{N+1}\times\R^{N+1}$ the following
estimates hold
\begin{itemize}
    \item [(i)]$\al|D_PW^\al(t,x,y,P)|\leq C_1$, where $C_1$ is introduced in (F2);
            \item [(ii)] $|\al W^\al(t,x,y,P)+\overline{F}(P)|\leq \al
    K_1$, where $K_1$ is introduced in Proposition \ref{barlesprop};
    \item [(iii)] $W^\al(t,x,y,0)\equiv0$;
    \item [(iv)] $\|D\overline{F}\|_\infty\leq C_1$.
    \end{itemize}
\end{lem}
\dim Let us fix $Q\in\R^{N+1}$. The Lipschitz continuity of $F$,
i.e. (F2), implies that the function $W(t,x,y)=W^\al(t,x,y,P+Q)$
satisfies
$$W_t+F(t,x,y,P+DW)+\al W\leq C_1|Q|$$ and then, by comparison
$$\al W(t,x,y)\leq \al W^\al(t,x,y,P)+C_1|Q|.$$
A similar argument shows that $\al W(t,x,y)\geq \al
W^\al(t,x,y,P)-C_1|Q|$. It then follows
$$\al |W^\al(t,x,y,P+Q)-W^\al(t,x,y,P)|\leq C_1|Q|,$$ which proves (i).

Let us turn out to (ii). We claim that
$$\mu:=\al \max_{\R\times\R^{N+1}} W^\al\geq -\overline{F}(P).$$ Indeed, $W^\al(t,x,y,P)$
is a supersolution of
$$W^\al_t+F(t,x,y,P+DW^\al)=-\mu.$$ Let $V$ be a bounded
subsolution of \eqref{V}, then by comparison between  $W^\al+\mu
t$ and $V-\overline{F}(P)t$, we have
$$V(t,x,y)-W^\al(t,x,y)\leq
V(0,x,y)-W^\al(0,x,y)+t(\overline{F}(P)+\mu).$$ Since $V$ and
$W^\al$ are bounded, dividing by $t>0$ and letting $t$ tend to
$+\infty$, we obtain $\mu\geq -\overline{F}(P)$. Then from (ii) of
Proposition \ref{barlesprop}, for $(t,x,y)\in\R\times\R^{N+1}$,
\beqs\begin{split}\al W^\al(t,x,y,P)&\geq
\al\min_{\R\times\R^{N+1}} W^\al\geq \al\max_{\R\times\R^{N+1}}
W^\al-\al K_1\geq -\overline{F}(P)-\al K_1.\end{split}\eeqs A
similar argument shows that
$$\al W^\al(t,x,y,P)+\overline{F}(P)\leq \al K_1;$$ this
concludes the proof of (ii).

Property (iii) follows from $F(t,x,y,0,0)=0$ and the uniqueness of
the periodic solution of \eqref{Wal}.

Finally, (iv) is an immediate consequence of
$$\overline{F}(P)-\overline{F}(Q)\leq
2\al K_1+\al\|D_PW^\al\|_\infty|P-Q|$$ and of (i). \finedim

We conclude this section by recalling some properties of the
solutions $u^0$ and $u^\ep$.

\begin{prop}\label{regulU0}There exist  constants $C_T, L>0$ such
that for any $(t,x),(s,y)\in [0,T]\times\R^N$
\beq\label{uepbound}|u^\ep(t,x)|,\,|u^0(t,x)|\leq C_T,\eeq
 \beq\label{ueplip}\begin{split}|u^0(t,x)-u^0(s,y)|\leq
L(|t-s|+|x-y|).\end{split}\eeq Moreover, for any $t\in[0,T]$, the
Lipschitz constant of $u^0(t,\cdot)$ is the Lipschitz constant of
the initial datum $u_0$.
\end{prop}
\dim By comparison $$|u^\ep(t,x)-u_0(x)|\leq C_0 t$$ where
$C_0=\max_{x,y,|p|\leq|u_0|_{1,\infty}}|H(x,y,p)|$. This implies
\eqref{uepbound} for $u^\ep$. Similarly can be showed the same
estimate for $u^0$.

The Lipschitz continuity of $u^0$ follows from the comparison
principle for \eqref{u0}, see \cite{bc}, Theorem III.3.7 and
Remark III.3.8.
 \finedim

\subsection{Proof of the main result} This section is devoted
to the proof of Theorem \ref{main}. We are going to show that for
any $T>0$
$$\sup_{[0,T]\times\R^{N+1}}|U^\ep(t,x,y)-U^0(t,x,y)|\leq
Ce^T\ep^\frac{1}{3},$$where $C$ does not depend on $T$. Since
$U^\ep(t,x,y)=u^\ep(t,x)-y$ and $U^0(t,x,y)=u^0(t,x)-y$, this
estimate automatically gives \eqref{main1}.

Let us consider a function $\phi:\R\rightarrow\R$  with the
following properties \beq\label{varphiprop}\left\{%
\begin{array}{ll}
    \phi'(s)>0, & \hbox{for any } s\in\R, \\
    \displaystyle \lim_{s\rightarrow+\infty}\phi(s)=1,\,\quad\quad\quad & \displaystyle \lim_{s\rightarrow-\infty}\phi(s)=0, \\
    |\phi(s)-\chi(s)|,\,|\phi'(s)|\leq \frac{K_2}{1+s^2}, & \hbox{for any }s\in\R,  \\
\end{array}%
\right. \eeq  where we have denoted by $\chi(s)$ the heaviside
function defined by \beqs \chi(s)=\left\{%
\begin{array}{ll}
    1, & \hbox{for }s\geq 0, \\
    0, & \hbox{for }s<0.
\end{array}%
\right. \eeqs For $n\in\N$, $\ep,\,\delta>0$, let us define the
function \beqs \varphi_\ep^{n,\delta}(s):=\sum_{i=-n}^n
\ep\phi\left(\frac{s-i\ep}{\delta}\right)-\ep (n+1).\eeqs Then we
have:
\begin{lem}\label{varphilem}Assume \eqref{varphiprop}. Then for any
$s\in\R$, the limit $\lim_{n\rightarrow+\infty}\varphi_\ep^{n,\delta}(s)$ exists and the function $\varphi_\ep^\delta$:
$$\varphi_\ep^\delta(s):=\lim_{n\rightarrow+\infty}\varphi_\ep^{n,\delta}(s)$$
 is of class $C^1$ with $(\varphi_\ep^\delta)'(s)>0$ for any $s\in\R$. Moreover
\beq\label{varphilim}\lim_{\delta\rightarrow0^+}\varphi_\ep^\delta(s)=\left\{%
\begin{array}{ll}
(i-1)\ep+\phi(0)\ep, & \hbox{if } s=i\ep,\\
    i\ep, & \hbox{if }i\ep<s<(i+1)\ep. \\
\end{array}%
\right.\eeq
\end{lem}
See the Appendix for the proof of the lemma.

Let us define
$$\widetilde{U}^{\ep,\delta}(t,x,y):=\varphi^\delta_\ep(U^\ep(t,x,y)).$$
Since $F$ satisfies the "geometrical" assumption (F4), the
function $\widetilde{U}^{\ep,\delta}$  is  solution of
\beq\label{Vep0}\left\{%
\begin{array}{ll}
    \widetilde{U}^{\ep,\delta}_t+F\left(\frac{t}{\ep},\frac{x}{\ep},\frac{U^\ep+y}{\ep},D_x \widetilde{U}^{\ep,\delta},
    D_y \widetilde{U}^{\ep,\delta}\right)=0,
     & (t,x,y)\in(0,T)\times\R^{N+1}, \\
    \widetilde{U}^{\ep,\delta}(0,x,y)=\varphi_{\ep}^\delta(u_0(x)-y), & (x,y)\in\R^{N+1}.\\
\end{array}%
\right.\eeq
 By stability
of viscosity solutions, see e.g. \cite{cil}, the limit
$\widetilde{U}^{\ep}(t,x,y)$  of
$\widetilde{U}^{\ep,\delta}(t,x,y)$ as $\delta\rightarrow0^+$ is
 a discontinuous viscosity solution of \eqref{Vep0} with
initial datum $\varphi_\ep(u_0(x)-y)$, where
$\varphi_\ep(s)=\lim_{\delta\rightarrow0^+}\varphi^\delta_\ep(s).$
 This means that
$(\widetilde{U}^\ep)^*=\limsup^*_{\delta\rightarrow0^+}\widetilde{U}^{\ep,\delta}$
(resp.
$(\widetilde{U}^\ep)_*=\liminf_{*\,\delta\rightarrow0^+}\widetilde{U}^{\ep,\delta}$)
is a viscosity subsolution (resp. supersolution) of \eqref{Vep0},
 and $(\widetilde{U}^\ep)^*(0,x,y)\leq
(\varphi_{\ep})^*(u_0(x)-y)$ (resp.
$(\widetilde{U}^\ep)_*(0,x,y)\geq (\varphi_{\ep})_*(u_0(x)-y)$).
Moreover, by \eqref{varphilim}
$$\widetilde{U}^{\ep}(t,x,y)=\left\{%
\begin{array}{ll}
    i\ep, & \hbox{if } i\ep<U^\ep(t,x,y)<(i+1)\ep ,\\
    (i-1)\ep+\phi(0)\ep, & \hbox{if }(t,x,y)\in \text{Int}\{U^\ep=i\ep\} .\\ \\
\end{array}%
\right. $$ At the points $(t,x,y)\in \partial\{U^\ep=i\ep\}$, the value of
$\widetilde{U}^{\ep}$ depends on the lower semi-continuous or the
upper semi-continuous envelope that we consider in the definition of
discontinuous viscosity solution. In particular, since $U^\ep$ is
continuous, $\widetilde{U}^{\ep}$ has the following properties
 \beq\label{varphidel-s}|(\widetilde{U}^\ep)^*(t,x,y)-U^\ep(t,x,y)|,\,|(\widetilde{U}^\ep)_*(t,x,y)-U^\ep(t,x,y)|\leq \ep\quad\text{for any }(t,x,y)\in[0,T]\times\R^{N+1}\eeq
 and \beq\label{0grad} D\widetilde{U}^\ep(t,x,y)=0\quad
\text{if }U^\ep(t,x,y)\neq i\ep,\,i\in\Z.\eeq Condition
\eqref{0grad} implies that $\widetilde{U}^\ep$ is actually a
solution of
\beqs\left\{%
\begin{array}{ll}
    \widetilde{U}^\ep_t+F\left(\frac{t}{\ep},\frac{x}{\ep},\frac{y}{\ep},D_x \widetilde{U}^\ep,D_y \widetilde{U}^\ep\right)=0,
    & (t,x,y)\in(0,T)\times\R^{N+1}, \\
    \widetilde{U}^\ep(0,x,y)=\varphi_\ep(u_0(x)-y), & (x,y)\in\R^{N+1}.\\
\end{array}%
\right. \eeqs Indeed, when $i\ep <U^\ep(t,x,y)<(i+1)\ep$, for some
$i\in\Z$, the function $\widetilde{U}^\ep$ is constant in a
neighborhood of $(t,x,y)$. Then the result follows from the fact
that $F(t,x,y,0)=0$. On the other hand, when $U^\ep(t,x,y)=i\ep$,
by periodicity,
$F\left(\frac{t}{\ep},\frac{x}{\ep},\frac{U^\ep+y}{\ep},P\right)=F\left(\frac{t}{\ep},\frac{x}{\ep},\frac{y}{\ep},P\right)$.

In order to estimate $|U^\ep-U^0|$ it is convenient to estimate
$|\widetilde{U}^\ep-U^0|$; indeed, $\frac{U^\ep}{\ep}$ does not any longer appear in the equation satisfied by
$\widetilde{U}^\ep$.%the dependence on $\frac{U^\ep}{\ep}$ disappears.

Let us define $V^{\ep}(t,x,y)=e^{-t}\widetilde{U}^{\ep}(t,x,y)$
and $V^0(t,x,y)=e^{-t}U^0(t,x,y)$. The functions $V^{\ep}$ and
$V^0$ are respectively solutions of

\beq\label{Vep}\left\{%
\begin{array}{ll}
    V^{\ep}_t+V^{\ep}+F\left(\frac{t}{\ep},\frac{x}{\ep},\frac{y}{\ep},D_xV^\ep,D_yV^\ep\right)=0,
     & (t,x,y)\in(0,T)\times\R^{N+1}, \\
    V^{\ep}(0,x,y)=\varphi_\ep(u_0(x)-y), & (x,y)\in\R^{N+1},\\
\end{array}%
\right. \eeq and

\beq\label{V0}\left\{%
\begin{array}{ll}
    V^0_t+V^0+\overline{F}(D_x V^0,D_yV^0)=0, & (t,x,y)\in(0,T)\times\R^{N+1}, \\
    V^0(0,x,y)=u_0(x)-y, & (x,y)\in\R^{N+1}.\\
\end{array}%
\right. \eeq For alleviating the notations, let us denote a vector
of $\R^{N+1}$ by $X=(x,x_{N+1})$, where $x\in\R^N$ and
$x_{N+1}\in\R$. We first estimate from above the difference
$(V^{\ep})^*-V^0$: for this, let us introduce the
auxiliary function
\beq\begin{split}\label{phi}\Phi(t,X,s,Y)&=(V^{\ep})^*(t,X)-V^0(s,Y)-\ep
W^\al\left(\frac{t}{\ep},\frac{X}{\ep},\frac{X-Y}{\ep^\beta}\right)\\&-\frac{|X-Y|^2}{2\ep^\beta}
-\frac{|t-s|^2}{2\sigma}
-\frac{r}{2}|X|^2-\frac{\eta}{T-t},\end{split}\eeq where
$\al=\ep^\theta$, $\theta,\,\beta,\,\sigma,\,r,\,\eta\in(0,1)$
will be fix later on and $\beta$ and $\theta$ satisfy
\beq\label{thetabeta}0<\theta<1-\beta.\eeq

In view of \eqref{uepbound},  \eqref{varphidel-s}, (i) of
Proposition \ref{barlesprop} and \eqref{F4},
$$\Phi(t,X,s,Y)\leq
2C_T+\ep+|x_{N+1}-y_{N+1}|+\frac{\ep}{\al}C_1\frac{|X-Y|}{\ep^\beta}-\frac{|X-Y|^2}{2\ep^\beta}-\frac{r}{2}|X|^2$$
for all $(t,X),\,(s,Y)\in[0,T]\times\R^{N+1}$. Hence, $\Phi$
attains a global maximum at some point
$(\ts,\Xs,\s,\Ys)\in([0,T]\times\R^{N+1})^2$. Standard arguments
show that $\ts,\,\s<T$ for $\sigma$ small enough.\\

{\bf Claim 1:} There exists a constant $M_1>0$ independent of
$\ep$ such that $\frac{|\ts-\s|}{\sigma}\leq M_1(1+|\ys_{N+1}|)$.
\\
The inequality
$\Phi(\ts,\Xs,\ts,\Ys)\leq\Phi(\ts,\Xs,\s,\Ys)$ and Proposition
\eqref{regulU0} imply
\beqs\begin{split}\frac{|\ts-\s|^2}{2\sigma}&\leq
V^0(\ts,\Ys)-V^0(\s,\Ys)\leq
|e^{-\ts}-e^{-\s}||U^0(\ts,\Ys)|+e^{-\s}|U^0(\ts,\Ys)-U^0(\s,\Ys)|\\&
\leq |\ts-\s|(C_T+|\ys_{N+1}|)+L|\ts-\s|\end{split}\eeqs from which Claim 1 follows.\\

{\bf Claim 2:} There exists a constant $M_2>0$ independent of
$\ep$ and $T$, such that $\frac{|\Xs-\Ys|}{\ep^\beta}\leq M_2$.
\\
The inequality $\Phi(\ts,\Xs,\s,\Xs)\leq\Phi(\ts,\Xs,\s,\Ys)$
implies \beqs\begin{split}\frac{|\Xs-\Ys|^2}{\ep^\beta}&\leq
V^0(\s,\Xs)-V^0(\s,\Ys)+\ep
W^\al\left(\frac{\ts}{\ep},\frac{\Xs}{\ep},0\right) -\ep
W^\al\left(\frac{\ts}{\ep},\frac{\Xs}{\ep},\frac{\Xs-\Ys}{\ep^\beta}\right).\end{split}\eeqs
Using \eqref{ueplip}, (i) of Lemma \ref{wproplem} and
\eqref{thetabeta} we then infer

\beqs \frac{|\Xs-\Ys|^2}{\ep^\beta}\leq
(L+1)|\Xs-\Ys|+\frac{\ep}{\al}C_1\frac{|\Xs-\Ys|}{\ep^\beta}=(L+1)|\Xs-\Ys|+\ep^{1-\theta-\beta}C_1
|\Xs-\Ys|\leq M_2 |\Xs-\Ys|.\eeqs  This concludes the proof of Claim 2.\\

{\bf Claim 3:} There exists a constant $M_3>0$ independent of
$\ep$ such that $r|\Xs|^2\leq M_3$.
\\ The inequality
$\Phi(\ts,0,\s,0)\leq\Phi(\ts,\Xs,\s,\Ys)$  implies
\beqs\begin{split}\frac{r}{2}|\Xs|^2&\leq
(V^\ep)^*(\ts,\Xs)-V^0(\s,\Ys)+V^0(\s,0)-(V^\ep)^*(\ts,0)+\ep
W^\al\left(\frac{\ts}{\ep},0,0\right)- \ep
W^\al\left(\frac{\ts}{\ep},\frac{\Xs}{\ep},\frac{\Xs-\Ys}{\ep^\beta}\right).\end{split}\eeqs
Then, using \eqref{uepbound}, \eqref{varphidel-s}, Claims 1 and 2,
(iii) of Lemma \ref{wproplem}, (i) of Proposition \ref{barlesprop}
and \eqref{F4}, we deduce
\beqs\begin{split}\frac{r}{2}|\Xs|^2&\leq
e^{-\ts}[U^\ep(\ts,\Xs)-U^0(\s,\Ys)]+|e^{-\ts}-e^{-\s}||U^0(\s,\Ys)|+\ep\\&+V^0(\s,0)-(V^\ep)^*(\ts,0)-
\ep
W^\al\left(\frac{\ts}{\ep},\frac{\Xs}{\ep},\frac{\Xs-\Ys}{\ep^\beta}\right)\\&
\leq
4C_T+M_2\ep^\beta+|\ts-\s|(C_T+|\ys_{N+1}|)+2\ep+\frac{\ep}{\al}C_1\frac{|\Xs-\Ys|}{\ep^\beta}\\&\leq
C+2\sigma M_1|\ys_{N+1}|^2\leq C+2\sigma M_1|\Xs|^2 ,\end{split}\eeqs and Claim 3 follows by choosing $\sigma<\frac{r}{8M_1}$.\\

Now, suppose first that $\ts=0$, then
\beqs\begin{split}&(V^\ep)^*(t,X)-V^0(t,X)-\ep
W^\al\left(\frac{t}{\ep},\frac{X}{\ep},0\right)-\frac{\eta}{T-t}-\frac{r}{2}|X|^2
\\&\leq (\varphi_\ep)^*(u_0(\xs)-\xs_{N+1})-V^0(\s,\Ys)- \ep
W^\al\left(0,\frac{\Xs}{\ep},\frac{\Xs-\Ys}{\ep^\beta}\right)\end{split}\eeqs
for any $(t,X)\in[0,T]\times\R^{N+1}$, from which, using (i) of
Proposition \ref{barlesprop}, (iii) of Lemma \ref{wproplem},
\eqref{F4} and Claim 2, we deduce
\beqs\begin{split}(V^\ep)^*(t,X)-V^0(t,X)\leq
(\varphi_\ep)^*(u_0(\xs)-\xs_{N+1})-V^0(\s,\Ys)+\frac{\eta}{T-t}+\frac{r}{2}|X|^2+\ep^{1-\theta}C_1M_2.\end{split}\eeqs

Letting $\sigma,\,\eta$ and $r$ go to $0^+$ and using
\eqref{varphidel-s} and Claim 2 we obtain
\beqs\begin{split}(V^\ep)^*(t,X)-V^0(t,X)&\leq
(\varphi_\ep)^*(u_0(\xs)-\xs_{N+1})-(u_0(\ys)-\ys_{N+1})+C\ep^{1-\theta}\\&\leq
(\varphi_\ep)^*(u_0(\xs)-\xs_{N+1})-(u_0(\xs)-\xs_{N+1})+(L+1)|\Xs-\Ys|+C\ep^{1-\theta}\\&\leq
C(\ep^\beta+\ep^{1-\theta})+\ep,\end{split}\eeqs which implies
\beq\label{t=0}U^\ep(t,X)-U^0(t,X)\leq
Ce^t(\ep^\beta+\ep^{1-\theta}).\eeq The same estimate can be
showed if $\s=0$.

Next, let us consider the case $\ts,\s>0$.

{\bf Claim 4:} There exists a constant $C>0$ independent of $\ep$
and $T$ such that
$$\frac{\ts-\s}{\sigma}+\frac{\eta}{(T-\ts)^2}+(V^\ep)^*(\ts,\Xs)
+\overline{F}\left(\frac{\Xs-\Ys}{\ep^\beta}\right) \leq
C(\ep^{1-\theta-\beta}+\ep^\theta).$$

The function
\beq\label{psiwrttx}\begin{split}(t,X)&\rightarrow
(V^\ep)^*(t,X)-\ep
W^\al\left(\frac{t}{\ep},\frac{X}{\ep},\frac{X-\Ys}{\ep^\beta}\right)-\frac{|X-\Ys|^2}{2\ep^\beta}-\frac{r}{2}|X|^2
-\frac{|t-\s|^2}{2\sigma} -\frac{\eta}{T-t}\end{split}\eeq has a
maximum at $(\ts,\Xs)$. By adding to $\Phi$ a smooth function
vanishing with its first derivative at $(\ts,\Xs)$, we may assume
the maximum is strict.

Next, for $j>0$, let us introduce the function
\beqs\begin{split}\Psi_j(t,s,X,Y,Z):&=(V^\ep)^*(t,X)-\ep W^\al
\left(s,Y,\frac{Z-\Ys}{\ep^\beta}\right)-\frac{|X-\Ys|^2}{2\ep^\beta}-\frac{r}{2}|X|^2
-\frac{|t-\s|^2}{2\sigma}\\& -\frac{\eta}{T-t}-\frac{j}{2}(|t-\ep
s|^2+|X-Z|^2+|X-\ep Y|^2).\end{split}\eeqs

Let $P_j=(t_j,s_j,X_j,Y_j,Z_j)$ be a maximum point of $\Psi_j$ on
the set
$$A:=\overline{B}(\ts,1)\times \overline{B}\left(\frac{\ts}{\ep},1\right)\times \overline{B}(\Xs,1)
\times
\overline{B}\left(\frac{\Xs}{\ep},1\right)\times\overline{B}(\Xs,1).$$
Since $(\ts,\Xs)$ is a strict maximum point of \eqref{psiwrttx},
$t_j\rightarrow\ts$, $s_j\rightarrow\frac{\ts}{\ep}$,
$X_j,\,Z_j\rightarrow\Xs$ and $Y_j\rightarrow\frac{\Xs}{\ep}$ as
$j\rightarrow +\infty$. Then, for $j$ large enough, $P_j$ lies
in the interior of $A$. Moreover, standard arguments show that
\beq\label{convergencej} j|t_j-\ep s_j|^2,\quad j|X_j-Z_j|^2,\quad
j|X_j-\ep Y_j|^2\rightarrow0\quad\text{as
}j\rightarrow+\infty.\eeq Remark that this implies in addition
that \beq\label{convergencej2}2j|t_j-\ep s_j||X_j-\ep Y_j|\leq
j|t_j-\ep s_j|^2+j|X_j-\ep Y_j|^2\rightarrow0\quad\text{as
}j\rightarrow+\infty.\eeq  Since $(V^\ep)^*$ and $W^\al$ are
respectively viscosity subsolutions of \eqref{Vep} and
supersolution of \eqref{Wal}, we obtain

\beq\begin{split}\label{usubsol}&\frac{t_j-\s}{\sigma}+\frac{\eta}{(T-t_j)^2}+j(t_j-\ep
s_j)+(V^\ep)^*(t_j, X_j)\\&+
F\left(\frac{t_j}{\ep},\frac{X_j}{\ep},\frac{X_j-\Ys}{\ep^\beta}+rX_j+j(X_j-Z_j)+j(X_j-\ep
Y_j)\right)\leq 0\end{split}\eeq and
\beq\begin{split}\label{wsupers}& j(t_j-\ep s_j)+\al
W^\al\left(s_j,Y_j,\frac{Z^j-\Ys}{\ep^\beta}\right)
+F\left(s_j,Y_j,\frac{Z^j-\Ys}{\ep^\beta}+j(X_j-\ep
Y_j)\right)\geq0.\end{split}\eeq Subtracting \eqref{usubsol} and
\eqref{wsupers} and using the Lipschitz continuity of $F$,
assumption (F2), we get
\beq\begin{split}\label{usub-wsuper}&\frac{t_j-\s}{\sigma}+\frac{\eta}{(T-t_j)^2}+(V^\ep)^*(t_j,
X_j)-\al
W^\al\left(s_j,Y_j,\frac{Z^j-\Ys}{\ep^\beta}\right)\\&\leq
\frac{C_1}{\ep}\left(\left|t_j-\ep s_j\right|+\left|X_j-\ep
Y_j\right|\right)\left(\frac{|Z^j-\Ys|}{\ep^\beta}+j|X_j-\ep
Y_j|\right)+C_1\left(\frac{|X_j-Z_j|}{\ep^\beta}+r|X_j|+j|X_j-Z_j|\right).\end{split}\eeq
 Let us estimate $j|X_j-Z_j|$. From the inequality
$\Psi_j(t_j,s_j,X_j,Y_j,X_j)\leq \Psi_j(t_j,s_j,X_j,Y_j,Z_j)$ we
deduce that
$$\frac{j}{2}|X_j-Z_j|^2\leq \ep W^\al
\left(s_j,Y_j,\frac{X^j-\Ys}{\ep^\beta}\right)-\ep W^\al
\left(s_j,Y_j,\frac{Z^j-\Ys}{\ep^\beta}\right),$$ and using (i) of
Lemma \ref{wproplem} we get
$$\frac{j}{2}|X_j-Z_j|^2\leq
C_1\frac{\ep}{\al}\frac{|X_j-Z_j|}{\ep^\beta}=C_1
\ep^{1-\theta-\beta}|X_j-Z_j|.$$ Then
\beq\label{z-x}j|X_j-Z_j|\leq 2C_1 \ep^{1-\theta-\beta}.\eeq Then,
passing to the limsup as $j\rightarrow+\infty$ in
\eqref{usub-wsuper} and taking into account Claim 2,
\eqref{convergencej} and \eqref{convergencej2}, we obtain

\beq\label{stimauepw}\frac{\ts-\s}{\sigma}+\frac{\eta}{(T-\ts)^2}+(V^\ep)^*(\ts,\Xs)-\al
W^\al\left(\frac{\ts}{\ep},\frac{\Xs}{\ep},\frac{\Xs-\Ys}{\ep}\right)\leq
C(\ep^{1-\theta-\beta}+r|\Xs|).\eeq By Claim 3, $r|\Xs|\leq
r^\frac{1}{2}M_3^\frac{1}{2}$, hence choosing $r>0$ such that
$r^\frac{1}{2}M_3^\frac{1}{2}\leq \ep^{1-\theta-\beta}$, we have
$r|\Xs|\leq \ep^{1-\theta-\beta}$.

Finally, Claim 4 easily follows from \eqref{stimauepw}, Claim 2
and the following inequality
$$-\al W^\al\left(\frac{\ts}{\ep},\frac{\Xs}{\ep},\frac{\Xs-\Ys}{\ep^\beta}\right)\geq
\overline{F}\left(\frac{\Xs-\Ys}{\ep^\beta}\right)-\al K_1\geq
\overline{F}\left(\frac{\Xs-\Ys}{\ep^\beta}\right)-K_1\ep^\theta$$
which comes from (ii) of Lemma \ref{wproplem} .\\

{\bf Claim 5:} There exists a constant $C>0$ independent of $\ep$
and $T$ such that
$$\frac{\ts-\s}{\sigma}+V^0(\s,
\Ys) +\overline{F}\left(\frac{\Xs-\Ys}{\ep^\beta}\right) \geq
-C\ep^{1-\theta-\beta}.$$

The function

\beqs\begin{split}(s,Y)&\rightarrow \phi(s,Y):=V^0(s,Y)+\ep
W^\al\left(\frac{\ts}{\ep},\frac{\Xs}{\ep},\frac{\Xs-Y}{\ep^\beta}\right)+\frac{|\Xs-Y|^2}{2\ep^\beta}
+\frac{|\ts-s|^2}{2\sigma}
\end{split}\eeqs
has a minimum at $(\s,\Ys)$, consequently $(0,0)\in
D^-\phi(\s,\Ys)$. If we set \beqs
\widetilde{V}(s,Y):=V^0(s,Y)+\frac{|\Xs-Y|^2}{2\ep^\beta}
+\frac{|\ts-s|^2}{2\sigma},\quad \widetilde{W}(Y):=\ep
W^\al\left(\frac{\ts}{\ep},\frac{\Xs}{\ep},\frac{\Xs-Y}{\ep^\beta}\right),\eeqs
 by properties of semijets of Lipschitz functions, see e.g. Lemma 2.4 in \cite{ci}, there exists $Q\in\R^{N+1}$ such that
  $$(0,Q)\in D^- \widetilde{V}(\s,\Ys)=D^-V^0(\s,\Ys)-\left(\frac{\ts-\s}{\sigma},\frac{\Xs-\Ys}{\ep^\beta}\right)
  \quad -Q\in D^-\widetilde{W}(\Ys).$$
Since $V^0$ is a supersolution of \eqref{V0}, we have
\beq\label{v0super}\frac{\ts-\s}{\sigma}+V^0(\s, \Ys)
+\overline{F}\left(\frac{\Xs-\Ys}{\ep^\beta}+Q\right) \geq 0.\eeq
By (i) of Lemma \ref{wproplem}, $$\left|\ep
W^\al\left(\frac{\ts}{\ep},\frac{\Xs}{\ep},\frac{\Xs-Y}{\ep^\beta}\right)-\ep
W^\al\left(\frac{\ts}{\ep},\frac{\Xs}{\ep},\frac{\Xs-Z}{\ep^\beta}\right)\right|
\leq
\frac{\ep}{\al}C_1\frac{|Y-Z|}{\ep^\beta}=C_1\ep^{1-\theta-\beta}|Y-Z|,$$
from which we get the following estimate of $Q$:
\beq\label{qlip}|Q|\leq C_1\ep^{1-\theta-\beta}.\eeq Then, Claim 5
follows from \eqref{v0super} using estimate \eqref{qlip} and the
 Lipschitz continuity of $\overline{F}$ assured by (iv) of Lemma \ref{wproplem}.\\

Claims 4 and 5 imply \beqs (V^\ep)^*(\ts, \Xs)-V^0(\s, \Ys)\leq
C(\ep^{1-\theta-\beta}+\ep^{\theta}),\eeqs for some constant $C$
independent of $\ep$ and $T$. Since $(\ts,\Xs,\s,\Ys)$ is a
maximum point of $\Phi$, we have
\beqs\begin{split}(V^\ep)^*(t,X)-V^0(t,X)\leq
\Phi(\ts,\Xs,\s,\Ys)+\ep
W^\al\left(\frac{t}{\ep},\frac{X}{\ep},0\right)+\frac{r}{2}|X|^2+\frac{\eta}{T-t},\end{split}\eeqs
for all $(t,X)\in[0,T]\times\R^{N+1}$. Then, by (iii) of Lemma
\ref{wproplem} \beqs\begin{split}(V^\ep)^*(t,X)-V^0(t,X)&\leq
(V^\ep)^*(\ts,\Xs)-V^0(\s,\Ys)-\ep
W^\al\left(\frac{\ts}{\ep},\frac{\Xs}{\ep},\frac{\Xs-\Ys}{\ep^\beta}\right)
+\frac{r}{2}|X|^2+\frac{\eta}{T-t}\\& \leq
C(\ep^{1-\theta-\beta}+\ep^{\theta})+
\frac{\ep}{\al}C_1\frac{|\Xs-\Ys|}{\ep^\beta}
+\frac{r}{2}|X|^2+\frac{\eta}{T-t}\\&\leq
C(\ep^{1-\theta-\beta}+\ep^{\theta})+\frac{r}{2}|X|^2+\frac{\eta}{T-t},
\end{split}\eeqs for some positive constant $C$.
Hence, sending $r,\,\eta,\rightarrow0^+$ and taking into account
\eqref{varphidel-s}, we get
$$U^\ep(t,X)-U^0(t,X)\leq
Ce^t(\ep^{1-\theta-\beta}+\ep^\theta).$$ Then, from the previous
estimate and \eqref{t=0}, we can conclude that for all
$\beta,\,\theta\in (0,1)$ satisfying \eqref{thetabeta} we have
$$U^\ep(t,X)-U^0(t,X)\leq
Ce^t(\ep^{1-\theta-\beta}+\ep^\theta+\ep^\beta),$$ for all
$(t,X)\in[0,T]\times\R^{N+1}$. The optimal choice of the
parameters is $\theta=\beta=\frac{1}{3}$, which gives
$$\sup_{[0,T]\times\R^{N+1}}(U^\ep(t,X)-U^0(t,X))\leq
C\ep^\frac{1}{3}.$$ The opposite inequality follows by similar
arguments, replacing $(V^\ep)^*$ with $V^0$ and $V^0$ with
$(V^\ep)_*$ in \eqref{phi}, and the proof of Theorem~\ref{main} in the
general case is complete.

Now, let us consider the case when $u_0$ is affine. Let us suppose that
$u_0(x)=p\cdot x+c_0$ for some $p\in\R^N$ and $c_0\in\R$. In this
case, the solution of \eqref{u0} is $u^0(t,x)=p\cdot
x+c_0-\overline{H}(p)t$. Let $\overline{V}$ be a bounded viscosity
supersolution of \eqref{V} with $p_x=p$ and $p_y=-1$.  Let us
define
$$V^\ep(t,X)= U^0(t,X)+\ep \overline{V}\left(\frac{t}{\ep},\frac{X}{\ep}\right). $$
Since $u_0(x)-y\geq\varphi_\ep(u_0(x)-y)-\ep$ then $V^\ep(0,X)\geq
\varphi_\ep(u_0(x)-y)-(M+1)\ep$ where $M=\|\overline{V}\|_\infty$.
Hence, it is easy to check that $V^\ep$ is a supersolution of
\beqs\left\{%
\begin{array}{ll}
    V^{\ep}_t+F\left(\frac{t}{\ep},\frac{X}{\ep},D_XV^\ep\right)=0,
     & (t,X)\in(0,T)\times\R^{N+1}, \\
    V^{\ep}(0,X)=\varphi_\ep(u_0(x)-y)-(M+1)\ep, & (x,y)\in\R^{N+1}.\\
\end{array}%
\right. \eeqs  By comparison we get $V^\ep(t,X)\geq
(\widetilde{U}^\ep)^*(t,X)-(M+1)\ep$ and this implies that
$U^0(t,X)-U^\ep(t,X)\geq -C\ep$. A similar argument shows that
$U^0(t,X)-U^\ep(t,X)\leq C\ep$ and this concludes the proof of the
theorem.\finedim

\section{Approximation of the effective Hamiltonian by Eulerian schemes}\label{sec:appr-effect-hamilt}
In this section we give an approximation of the effective
Hamiltonian $\overline{F}(P)$.
 To this end, we introduce an approximation scheme for the equation \eqref{Wal} and for simplicity we only discuss  the case $N=2$.
Given $N_X$ and $N_t$ positive integers, we introduce $\Delta t=1/N_t$,
$h=1/N_X$ and
$$\R^2_h:=\{X_{i,j}=(x_i,y_j)\,|\,x_i=ih,\,y_j=jh,\,i,j\in\Z\},$$
$$\R_{\Delta t}:=\{t_n=n\Delta t\,|\,n\in\Z\}.$$
An anisotropic mesh with steps $h_1$ and $h_2$  is possible too; we take $h_1=h_2$ only for simplicity. We denote by
$W_{i,j}^{n,P,\al}$ our numerical approximation of $W^{P,\al}$ at
$(t_n,x_i,y_j)\in\R_{\Delta t}\times\R^2_h$. For \eqref{Wal} we
consider the implicit Eulerian scheme of the form

\beq\label{eulsch}
\frac{W_{i,j}^{n+1,P,\al}-W_{i,j}^{n,P,\al}}{\Delta t}+\al
W_{i,j}^{n+1,P,\al}+S(t_n,x_i,y_j,h,[W^{n+1,P,\al}]_{i,j})=0,\eeq
where \beq\label{S}\begin{split}
S&(t_n,x_i,y_j,h,[W]_{i,j})\\&=g(t_n,x_i,y_j,(\Delta_1^+W)_{i,j}+p_x,(\Delta_1^+W)_{i-1,j}+p_x,(\Delta_2^+W)_{i,j}+p_y,
(\Delta_2^+W)_{i,j-1}+p_y)\end{split}\eeq and \beqs
(\Delta_1^+W)_{i,j}=\frac{W_{i+1,j}-W_{i,j}}{h},\quad
(\Delta_2^+W)_{i,j}=\frac{W_{i,j+1}-W_{i,j}}{h}.\eeqs We make the
following assumptions on $g$:
\renewcommand{\labelenumi}{(g\arabic{enumi})}
\begin{enumerate}
    \item Monotonicity: $g$ is nonincreasing with respect to its
    fourth and sixth arguments, and  nondecreasing with respect to its
    fifth and seventh arguments;
    \item Consistency: for any $t\in\R$, $(x,y)\in\R^2$ and $(q_x,q_y)\in\R^2$
     $$g(t,x,y,q_x,q_x,q_y,q_y)=F(t,x,y,q_x,q_y).$$
    \item Periodicity: for any $t\in\R$, $(x,y)\in\R^2$ and $Q\in\R^4$
    $$g(t+1,x+1,y+1,Q)=g(t,x,y,Q);$$
    \item Regularity: $g$ is locally Lipschitz continuous and there exists $\widetilde{C}_1>0$
     such that for any $t\in\R$, $(x,y)\in\R^2$ and $Q\in\R^4$
    $$|D_{Q}g(t,x,y,Q)|\leq \widetilde{C}_1;$$
    \item Coercivity: there exist $\widetilde{C}_2,\,\widetilde{C}_3>0$ such that for any $t\in\R$, $(x,y)\in\R^2$, $(q_1,q_2)\in\R^2$
$$g(t,x,y,q_1,q_2,0,0)\geq \widetilde{C}_2(|q_1^-|^2+|q_2^+|^2)^\frac{1}{2}-\widetilde{C}_3;$$
\item For any $t\in\R$, $(x,y_1),(x,y_2)\in\R^2$,  $q_1,q_2\in\R$
$$g(t,x,y_1,q_1,q_2,0,0)=g(t,x,y_2,q_1,q_2,0,0).$$
\end{enumerate}
The points (g1)-(g4) are standard assumptions in the study of
numerical schemes for Hamilton-Jacobi equations. The coercivity
hypothesis  (g5) can be substituted by the weaker condition \beqs
\lim_{q_1^++q_2^-\rightarrow+\infty}g(x,y,q_1,q_2,q_3,q_4)=+\infty\eeqs
if $g$ (and hence $F$) does not depend on time. If $g$ is
homogeneous of degree 1 w.r.t. $Q$, then the two coercivity
conditions are equivalent.

 As an example, we suppose that the Hamiltonian $F$ is of the form $F(t,x,y,p_x,p_y)=a(t,x)|p_x|+b(t,x,y)|p_y|$, with $a$ and $b$
 Lipschitz continuous functions and $a(t,x)\geq \widetilde{C}_2>0$; we consider a generalization of the Godunov scheme proposed in \cite{os}:
\beqs\begin{split}&g(t,x,y,q_1,q_2,q_3,q_4)\\&=a(t,x)[(q_1^-)^2+(q_2^+)^2]^\frac{1}{2}+b^+(t,x,y)[(q_3^-)^2+(q_4^+)^2]^\frac{1}{2}
-b^-(t,x,y)[(q_3^+)^2+(q_4^-)^2]^\frac{1}{2}.\end{split}\eeqs
%[(\min(q_1,0))^2+(\max(q_2,0))^2]^\frac{1}{2}+b(t,x,y)
where $q^+=\max(q,0)$ and $q^-=(-q)^+$. Then hypothesis (g1)-(g6)
are satisfied.

The following theorem is the discrete version of the analogous
result in \cite{b} for the exact solution $W^{P,\al}$ of \eqref{Wal}.
\begin{thm}\label{ergodicapprox}Assume (g1)-(g6). Then we have
\begin{itemize}
\item [(i)]For any $P=(p_x,p_y)\in\R^2$, $\al,\,h,\,\Delta t>0$ there
exists a unique $(W_{i,j}^{n,P,\al})$ periodic solution of
\eqref{eulsch};
\item [(ii)] There exists a constant $\widetilde{K}_1$ depending on
$\|F(\cdot,\cdot,\cdot,P)\|_\infty$, $\widetilde{C}_1$ in (g4),
$\widetilde{C}_2,\,\widetilde{C}_3$ in (g5), $p_x$ and $p_y$, but
independent of $\al$, $h$ and $\Delta t$
 such that
\beqs\max_{i,j,n}W_{i,j}^{n,P,\al}-\min_{i,j,n}W_{i,j}^{n,P,\al}\leq
\widetilde{K}_1;\eeqs
\item [(iii)] There exists a constant
$\overline{F}_h^{\Delta t}(P)$ such that
\beq\label{approxeffham}\lim_{\al\rightarrow0^+} \al
W_{i,j}^{n,P,\al}=-\overline{F}_{h}^{\Delta t}(P)\quad \forall
i,j,n;\eeq
\item [(iv)]$\overline{F}_{h}^{\Delta t}(P)$ is the unique
number $\lams_h^{\Delta t}\in\R$ such that the equation
\beq\label{ergodiceq} \frac{W_{i,j}^{n+1,P}-W_{i,j}^{n,P}}{\Delta
t}+S(t_n,x_i,y_j,h,[W^{n+1,P}]_{i,j})=\lams_h^{\Delta t}\eeq
admits a bounded solution.
\end{itemize}
\end{thm}
\dim A proof of the existence of a unique solution of
\eqref{eulsch} in the uniform grid on the torus with step $h$ is
given in \cite{cq}.

Let us prove (ii). First, remark that by comparison with constants
we have \beq\label{wc0}\max_{i,j,n}|\al W_{i,j}^{n,P,\al}|\leq
C_0,\eeq where $C_0:=\|F(\cdot,\cdot,\cdot,P)\|_\infty.$ Next, let
us define
$$\overline{W}_i^{n}:=\max_{j}W_{i,j}^{n,P,\al}.$$
We claim that $\overline{W}_i^{n}$ satisfies \beqs
\frac{\overline{W}_{i}^{n+1}-\overline{W}_{i}^{n}}{\Delta t}+\al
\overline{W}_i^{n+1}+\overline{S}(t_n,x_i,h,[\overline{W}^{n+1}]_i)\leq
0,\eeqs where
$$\overline{S}(t_n,x_i,h,[W]_i):=\min_{j}
g(t_n,x_i,y_j,(\Delta_1^+W)_i+p_x,(\Delta_1^+W)_{i-1}+p_x,p_y,p_y).$$
Indeed, for any $i$ and $n$, denote by $\overline{j}_{(i,n)}$ the
index $j$ such that
$\overline{W}_i^{n}=\max_{j}W_{i,j}^{n,P,\al}=W_{i,\overline{j}_{(i,n)}}^{n,P,\al}$,
then

\beqs
\frac{W_{i,\overline{j}_{(i,n+1)}}^{n+1,P,\al}-W_{i,\overline{j}_{(i,n+1)}}^{n,P,\al}}{\Delta
t}\geq
\frac{W_{i,\overline{j}_{(i,n+1)}}^{n+1,P,\al}-W_{i,\overline{j}_{(i,n)}}^{n,P,\al}}{\Delta
t}=\frac{\overline{W}_{i}^{n+1}-\overline{W}_{i}^{n}}{\Delta
t},\eeqs

\beqs
(\Delta_1^+W^{n+1,P,\al})_{i,\overline{j}_{(i,n+1)}}=\frac{W_{i+1,\overline{j}_{(i,n+1)}}^{n+1,P,\al}-W_{i,\overline{j}_{(i,n+1)}}^{n+1,P,\al}}{h}\leq
\frac{W_{i+1,\overline{j}_{i+1,n+1}}^{n+1,P,\al}-W_{i,\overline{j}_{(i,n+1)}}^{n+1,P,\al}}{h}=
(\Delta_1^+\overline{W}^{n+1})_i,\eeqs

 \beqs
(\Delta_1^+W^{n+1,P,\al})_{i-1,\overline{j}_{(i,n+1)}}=\frac{W_{i,\overline{j}_{(i,n+1)}}^{n+1,P,\al}-W_{i-1,\overline{j}_{(i,n+1)}}^{n+1,P,\al}}{h}\geq
\frac{W_{i,\overline{j}_{(i,n+1)}}^{n+1,P,\al}-W_{i-1,\overline{j}_{(i-1,n+1)}}^{n+1,P,\al}}{h}=
(\Delta_1^+\overline{W}^{n+1})_{i-1},\eeqs and \beqs (\Delta_2^+
W^{n+1,P,\al})_{i,\overline{j}_{(i,n+1)}}=\frac{W_{i,\overline{j}_{(i,n+1)}+1}^{n+1,P,\al}-W_{i,\overline{j}_{(i,n+1)}}^{n+1,P,\al}}{h}\leq
0,\eeqs \beqs (\Delta_2^+
W^{n+1,P,\al})_{i,\overline{j}_{(i,n+1)}-1}=\frac{W_{i,\overline{j}_{(i,n+1)}}^{n+1,P,\al}-W_{i,\overline{j}_{(i,n+1)}-1}^{n+1,P,\al}}{h}\geq
0.\eeqs Since $(W_{i,j}^{n,P,\al})$ satisfies \eqref{eulsch},
using the monotonicity assumption (g1), we get \beqs\begin{split}
&\frac{\overline{W}_{i}^{n+1}-\overline{W}_{i}^{n}}{\Delta t}+\al
\overline{W}_i^{n+1}+\overline{S}(t_n,x_i,h,[\overline{W}^{n+1}]_i)\\&
\leq \frac{\overline{W}_{i}^{n+1}-\overline{W}_{i}^{n}}{\Delta
t}+\al
W_{i,\overline{j}_{(i,n+1)}}^{n+1,P,\al}\\&+g(t_n,x_i,y_{\overline{j}_{(i,n+1)}},(\Delta_1^+\overline{W}^{n+1})_i
+p_x,(\Delta_1^+\overline{W}^{n+1})_{i-1}+p_x,p_y,p_y)
\\&\leq \frac{W_{i,\overline{j}_{(i,n+1)}}^{n+1,P,\al}-W_{i,\overline{j}_{(i,n+1)}}^{n,P,\al}}{\Delta
t}+ \al W_{i,\overline{j}_{(i,n+1)}}^{n+1,P,\al}\\&+
g(t_n,x_i,y_{\overline{j}_{(i,n+1)}},(\Delta_1^+W^{n+1,P,\al})_{i,\overline{j}_{(i,n+1)}}+p_x,
(\Delta_1^+W^{n+1,P,\al})_{i-1,\overline{j}_{(i,n+1)}}+p_x,\\&(\Delta_2^+
W^{n+1,P,\al})_{i,\overline{j}_{(i,n+1)}}+p_y,(\Delta_2^+
W^{n+1,P,\al})_{i,\overline{j}_{(i,n+1)}-1}+p_y)\\& \leq
0,\end{split}\eeqs as desired. Then, by (g4), (g5) and
\eqref{wc0}, we see that $\overline{W}_{i}^{n}$ satisfies
\beqs\frac{\overline{W}_{i}^{n+1}-\overline{W}_{i}^{n}}{\Delta
t}+\widetilde{C}_2\left(|[(\Delta_1^+\overline{W}^{n+1})_i
+p_x]^-|^2+|[(\Delta_1^+\overline{W}^{n+1})_{i-1}+p_x]^+|^2\right)^\frac{1}{2}-\leq
0,\eeqs where $K_1=C_0+\widetilde{C}_3+2\widetilde{C}_1|p_y|$. In
particular we infer that
$$\overline{W}_{i}^{n+1}-\overline{W}_{i}^{n}\leq K_1\Delta
t,$$ which implies that if $n\geq m$ then
\beq\label{wtempo}\overline{W}_{i}^{n}-\overline{W}_{i}^{m}\leq
K_1(n-m)\Delta t=K_1(t_n-t_m).\eeq

Next, let us consider $$\overline{\overline{W}}_i=\max_n
\overline{W}_i^n.$$ Similar arguments as before show that
$\overline{\overline{W}}_i$ satisfies

$$\widetilde{C}_2\left(|[(\Delta_1^+\overline{\overline{W}})_i
+p_x]^-|^2+|[(\Delta_1^+\overline{\overline{W}})_{i-1}+p_x]^+|^2\right)^\frac{1}{2}\leq
K_1,$$

which implies the existence of a constant $K_2>0$ depending on
$C_0,\widetilde{C}_1,\widetilde{C}_2,\widetilde{C}_3,p_x$ and
$p_y$ such that
\beq\label{liwoveover}\max_i|(\Delta_1^+\overline{\overline{W}})_i|\leq
K_2.\eeq

Now, let $(i_1,n_1)$ and $(i_2,n_2)$ be such that
$\max_{i,n}\overline{W}_i^n=\overline{W}_{i_1}^{n_1}$ and
$\min_{i,n}\overline{W}_i^n=\overline{W}_{i_2}^{n_2}$, and let
$n_{i_2}$ be such that $\overline{\overline{W}}_{i_2}=\max_n
\overline{W}_{i_2}^n=\overline{W}_{i_2}^{n_{i_2}}$. By
periodicity, we may take $|x_{i_1}-x_{i_2}|\leq 1$ and $0\leq
t_{n_{i_2}}-t_{n_2}\leq 1$. Then using \eqref{liwoveover} and
\eqref{wtempo}, we get
\beqs\begin{split}\overline{W}_{i_1}^{n_1}&=\overline{\overline{W}}_{i_1}\\&\leq
\overline{\overline{W}}_{i_2}+K_2|x_{i_1}-x_{i_2}|\\&\leq\overline{W}_{i_2}^{n_{i_2}}+K_2\\&\leq
\overline{W}_{i_2}^{n_{2}}+K_1(t_{n_{i_2}}-t_{n_2})+K_2\\&\leq
\overline{W}_{i_2}^{n_{2}}+K_0.\end{split}\eeqs Then we have
proved that
\beq\label{maxminoverw}\max_{i,n}\overline{W}_i^n-\min_{i,n}\overline{W}_i^n\leq
K_0,\eeq where $K_0$ depends only on
$C_0,\widetilde{C}_1,\widetilde{C}_2,\widetilde{C}_3,p_x$ and
$p_y$.

 Next, we consider the behavior of
$W_{i,j}^{n,P,\al}$ in $j$. We claim that \beqs
W_{i,j_1}^{n,P,\al}+p_yy_{j_1}\leq
W_{i,j_2}^{n,P,\al}+p_yy_{j_2}\quad \text{if }j_1\geq j_2\text{
and }p_y<0,\eeqs \beq\label{monp<0} W_{i,j_1}^{n,P,\al}=
W_{i,j_2}^{n,P,\al}\quad \text{for any }j_1,j_2\text{ if
}p_y=0,\eeq \beqs W_{i,j_1}^{n,P,\al}+p_yy_{j_1}\geq
W_{i,j_2}^{n,P,\al}+p_yy_{j_2}\quad \text{if }j_1\geq j_2\text{
and }p_y>0.\eeqs Let us consider the case $p_y<0$. Suppose by
contradiction that
$$M:=\max_{i,n,j_1\geq
j_2}(W_{i,j_1}^{n,P,\al}-W_{i,j_2}^{n,P,\al}+p_y(y_{j_1}-y_{j_2}))=W_{\overline{i},\overline{j}_1}^{\overline{n},P,\al}
-W_{\overline{i},\overline{j}_2}^{\overline{n},P,\al}+p_y(y_{\overline{j}_1}-y_{\overline{j}_2})>0.$$
Then $\overline{j}_1\geq \overline{j}_2+1$. We have the following
estimate \beqs\begin{split} (\Delta_1^+
W^{\overline{n},P,\al})_{\overline{i},\overline{j}_1}-(\Delta_1^+
W^{\overline{n},P,\al})_{\overline{i},\overline{j}_2}&=\frac{W^{\overline{n},P,\al}_{\overline{i}+1,\overline{j}_1}
-W^{\overline{n},P,\al}_{\overline{i},\overline{j}_1}}{h}-
\frac{W^{\overline{n},P,\al}_{\overline{i}+1,\overline{j}_2}-W^{\overline{n},P,\al}_{\overline{i},\overline{j}_2}}{h}\\&=
\frac{W^{\overline{n},P,\al}_{\overline{i}+1,\overline{j}_1}-W^{\overline{n},P,\al}_{\overline{i}+1,\overline{j}_2}}{h}
-\frac{W^{\overline{n},P,\al}_{\overline{i},\overline{j}_1}-W^{\overline{n},P,\al}_{\overline{i},\overline{j}_2}}{h}\leq
0.\end{split}\eeqs Similarly
$$(\Delta_1^+
W^{\overline{n},P,\al})_{\overline{i}-1,\overline{j}_1}\geq(\Delta_1^+
W^{\overline{n},P,\al})_{\overline{i}-1,\overline{j}_2},$$ and
$$\frac{W^{\overline{n},P,\al}_{\overline{i},\overline{j}_1}
-W^{\overline{n}-1,P,\al}_{\overline{i},\overline{j}_1}}{\Delta
t}\geq
\frac{W^{\overline{n},P,\al}_{\overline{i},\overline{j}_2}-W^{\overline{n}-1,P,\al}_{\overline{i},\overline{j}_2}}{\Delta
t}.$$

Moreover, we have \beqs\begin{split}
(\Delta_2^+W^{\overline{n},P,\al})_{\overline{i},\overline{j}_1}+p_y&=\frac{W^{\overline{n},P,\al}_{\overline{i},\overline{j}_1+1}
-W^{\overline{n}P,\al}_{\overline{i},\overline{j}_1}}{h}+p_y
\\&=\frac{W^{\overline{n},P,\al}_{\overline{i},\overline{j}_1+1}-W^{\overline{n},P,\al}_{\overline{i},\overline{j}_2}}{h}
+p_y \frac {y_{\overline{j}_1+1}-y_{\overline{j}_2}} h
-\frac{W^{\overline{n},P,\al}_{\overline{i},\overline{j}_1}-W^{\overline{n},P,\al}_{\overline{i},\overline{j}_2}}{h}
-p_y\frac {y_{\overline{j}_1}-y_{\overline{j}_2}} h \leq 0,\end{split}\eeqs
similarly \beqs
(\Delta_2^+W^{\overline{n},P,\al})_{\overline{i},\overline{j}_1-1}+p_y\geq
0,\quad
(\Delta_2^+W^{\overline{n},P,\al})_{\overline{i},\overline{j}_2}+p_y\geq
0,\quad
(\Delta_2^+W^{\overline{n},P,\al})_{\overline{i},\overline{j}_2-1}+p_y\leq
0.\eeqs Then, since $W_{i,j}^{\overline{n},P,\al}$ satisfies
\eqref{eulsch}, using assumptions (g1) and (g6), we get \beqs
\begin{split}\al (W_{\overline{i},\overline{j}_1}^{\overline{n},P,\al}
-W_{\overline{i},\overline{j}_2}^{\overline{n},P,\al})&\leq
-g(t_{\overline{n}},x_{\overline{i}},y_{\overline{j}_1},(\Delta_1^+
W^{\overline{n},P,\al})_{\overline{i},\overline{j}_1}+p_x,(\Delta_1^+
W^{\overline{n},P,\al})_{\overline{i}-1,\overline{j}_1}+p_x,
0,0)\\&+g(t_{\overline{n}},x_{\overline{i}},y_{\overline{j}_2},(\Delta_1^+
W^{\overline{n},P,\al})_{\overline{i},\overline{j}_1}+p_x,(\Delta_1^+
W^{\overline{n},P,\al})_{\overline{i}-1,\overline{j}_1}+p_x,
0,0)=0.\end{split}\eeqs This implies that $$0<\al
M=\al(W_{\overline{i},\overline{j}_1}^{\overline{n},P,\al}
-W_{\overline{i},\overline{j}_2}^{\overline{n},P,\al}+p_y(y_{\overline{j}_1}-y_{\overline{j}_2}))\leq
\al p_y(y_{\overline{j}_1}-y_{\overline{j}_2})<0,$$ which is a
contradiction and this concludes the proof of \eqref{monp<0} for
$p_y<0$. The case $p_y\geq 0$ can be treated in an analogous way.

Now, to prove (ii), we use the properties \eqref{maxminoverw} and
\eqref{monp<0} of $W_{i,j}^{n,P,\al}$ and again we only consider
the case $p_y<0$. Let $(i_1,j_1,n_1)$ and $(i_2,j_2,n_2)$ be such
that $W_{i_1,j_1}^{n_1,P,\al}=\max_{i,j,n}W_{i,j}^{n,P,\al}$ and
$W_{i_2,j_2}^{n_2,P,\al}=\min_{i,j,n}W_{i,j}^{n,P,\al}$. Let
$\overline{j}$ be such that
$\overline{W}_{i_2}^{n_2}=W_{i_2,\overline{j}}^{n_2,P,\al}$. By
periodicity, we can take $0\leq y_{\overline{j}}-y_{j_2}\leq 1$
and $|x_{i_1}-x_{i_2}|\leq1$. Then \beqs\begin{split}
W_{i_1,j_1}^{n_1,P,\al}&=\overline{W}_{i_1}^{n_1}\\&\leq
\overline{W}_{i_2}^{n_2}+K_0\\&=W_{i_2,\overline{j}}^{n_2,P,\al}+K_0\\&\leq
W_{i_2,j_2}^{n_2,P,\al}+p_y(y_{j_2}-y_{\overline{j}})+K_0\\&\leq
W_{i_2,j_2}^{n_2,P,\al}-p_y+K_0,\end{split}\eeqs and this
concludes the proof of (ii).

The property (iii) easily follows from (ii) and \eqref{wc0}.
Indeed, from \eqref{wc0}, up to subsequence,
$\al\min_{i,j,n}W_{i,j}^{n,P,\al}$ converges to a constant
$-\overline{F}_{h}^{\Delta t}(P)$ as $\al\rightarrow0^+$. Then
from (ii), for any $i,j,n$, we get \beqs\begin{split}|\al
W_{i,j}^{n,P,\al}+\overline{F}_{h}^{\Delta t}(P)|&\leq
|\al\min_{i,j,n}W_{i,j}^{n,P,\al}+\overline{F}_{h}^{\Delta
t}(P)|+\al|W_{i,j}^{n,P,\al}-\min_{i,j,n}W_{i,j}^{n,P,\al}|\\&\leq
|\al\min_{i,j,n}W_{i,j}^{n,P,\al}+\overline{F}_{h}^{\Delta t}|+\al
\widetilde{K}_1\rightarrow 0\quad\text{as
}\al\rightarrow0^+,\end{split}\eeqs and (iii) is proved.

Let us turn to (iv). Let us define
$Z_{i,j}^{n,P,\al}=W_{i,j}^{n,P,\al}-\min_{i,j,n}W_{i,j}^{n,P,\al}$.
By (ii), up to subsequence, $(Z_{i,j}^{n,P,\al})$ converges to a
grid function $(Z_{i,j}^{n,P})$ as $\al\rightarrow0^+$. The grid
function $(Z_{i,j}^{n,P,\al})$ satisfies \beqs
\frac{Z_{i,j}^{n+1,P,\al}-Z_{i,j}^{n,P,\al}}{\Delta t}+\al
Z_{i,j}^{n+1,P,\al}+S(t_n,x_i,y_j,h,[Z^{n+1,P,\al}]_{i,j})=-\al
\min_{i,j,n}W_{i,j}^{n,P,\al}.\eeqs Letting $\al\rightarrow0^+$,
since by (ii) $(Z_{i,j}^{n,P,\al})$ is bounded and
$\al\min_{i,j,n}W_{i,j}^{n,P,\al}\rightarrow-\overline{F}_h^{\Delta
t}$, we see that $(Z_{i,j}^{n,P})$ is a solution of
\eqref{ergodiceq} with $\lams_h^{\Delta t}=\overline{F}_h^{\Delta
t}$.

To prove the uniqueness of a  solution $(\lams_h^{\Delta
t},(W_{i,j}^{n,P}))$   of \eqref{ergodiceq}, we show that
if there exists a subsolution $(U_{i,j}^{n,P})$ of
\eqref{ergodiceq} with $\lams_h^{\Delta t}=\lam_1$ and a
supersolution $(V_{i,j}^{n,P})$ of \eqref{ergodiceq} with
$\lams_h^{\Delta t}=\lam_2$, then $\lam_2\leq\lam_1$.

Let
$M=\max_{i,j,n}(U_{i,j}^{n,P}-V_{i,j}^{n,P})=U_{i_0,j_0}^{n_0,P}-V_{i_0,j_0}^{n_0,P}$.
Then
$$\frac{U_{i_0,j_0}^{n_0,P}-U_{i_0,j_0}^{n_0-1,P}}{\Delta t}\geq \frac{V_{i_0,j_0}^{n_0,P}-V_{i_0,j_0}^{n_0-1,P}}{\Delta
t},$$
$$(\Delta_1^+U^{n_0,P})_{i_0,j_0}\leq
(\Delta_1^+V^{n_0,P})_{i_0,j_0},\quad
(\Delta_1^+U^{n_0,P})_{i_0-1,j_0}\geq
(\Delta_1^+V^{n_0,P})_{i_0-1,j_0},$$
$$(\Delta_2^+U^{n_0,P})_{i_0,j_0}\leq
(\Delta_2^+V^{n_0,P})_{i_0,j_0},\quad
(\Delta_2^+U^{n_0,P})_{i_0,j_0-1}\geq
(\Delta_2^+V^{n_0,P})_{i_0,j_0-1}.$$ From the monotonicity of $g$,
\beqs\begin{split}\lam_1&\geq
\frac{U_{i_0,j_0}^{n_0,P}-U_{i_0,j_0}^{n_0-1,P}}{\Delta t}+
g\left(t_{n_0},x_{i_0},y_{j_0},(\Delta_1^+U^{n_0,P})_{i_0,j_0}+p_x,(\Delta_1^+U^{n_0,P})_{i_0-1,j_0}+p_x,\right.\\&
\left.(\Delta_2^+U^{n_0,P})_{i_0,j_0}+p_y,
(\Delta_2^+U^{n_0,P})_{i_0,j_0-1}+p_y\right)\\& \geq
\frac{V_{i_0,j_0}^{n_0,P}-V_{i_0,j_0}^{n_0-1,P}}{\Delta t} +
g\left(t_{n_0},x_{i_0},y_{j_0},(\Delta_1^+V^{n_0,P})_{i_0,j_0}+p_x,(\Delta_1^+V^{n_0,P})_{i_0-1,j_0}+p_x,\right.\\&
\left.(\Delta_2^+V^{n_0,P})_{i_0,j_0}+p_y,
(\Delta_2^+V^{n_0,P})_{i_0,j_0-1}+p_y\right)\\&\geq\lam_2.\end{split}\eeqs
This concludes the proof of (iv).
 \finedim
We need a more precise estimate on the rate of convergence of $\al
W_{i,j}^{n,\al,P}$ to $\overline{F}_h^{\Delta t}(P)$:
 \begin{prop}\label{wal+fberapprox} Assume (g1)-(g6). Then for any
 $i,j,n$
 $$|\al W_{i,j}^{n,\al,P}+\overline{F}_h^{\Delta
t}(P)|\leq \widetilde{K}_1\al,$$ where
$\widetilde{K}_1=\widetilde{K}_1(P)$ is the constant in (ii) of
Theorem \ref{ergodicapprox}.
 \end{prop} \dim As in  the proof
of (ii) of Lemma \ref{wproplem}, the result follows from the
comparison principle for \eqref{eulsch} and (ii) of Theorem
\ref{ergodicapprox}.\finedim

Now, we are ready to show that the function
$\overline{F}_h^{\Delta t}$ is actually an approximation of the
effective Hamiltonian $\overline{F}$.
\begin{prop}\label{convappham} Assume
(g1)-(g6). Let $\overline{F}_h^{\Delta t}$ be defined by
\eqref{approxeffham} and let $\overline{F}$ be the effective
Hamiltonian. Then, for any $P\in\R^{2}$ \beqs\lim_{(\Delta t,
h)\rightarrow(0,0)}\overline{F}_h^{\Delta
t}(P)=\overline{F}(P)\eeqs uniformly on compact sets of $\R^{2}$.
\end{prop}
\dim To show the result we estimate
$W^{P,\al}(t_n,x_i,y_j)-W_{i,j}^{n,P,\al}$. To this end, following
the same proof as in \cite{cl} and \cite{acc}, we assume that
$$\sup_{i,j,n}|\al W^{P,\al}(t_n,x_i,y_j)-\al W_{i,j}^{n,P,\al}|=\sup_{i,j,n}(\al W^{P,\al}(t_n,x_i,y_j)
-\al W_{i,j}^{n,P,\al})=m\geq0.$$ The case when $\sup_{i,j,n}|\al
W^{P,\al}(t_n,x_i,y_j)-\al W_{i,j}^{n,P,\al}|=\sup_{i,j,n}(\al
W_{i,j}^{n,P,\al}-\al W^{P,\al}(t_n,x_i,y_j))$ is handled in a
similar manner.

 For simplicity of
notations we omit the index $P$. Let us denote $W_{h,\Delta
t}^{\al }(t_n,X_{i,j}):=W_{i,j}^{n,\al}$,
$(t_n,X_{i,j})\in\R_{\Delta t}\times\R^2_h$. For
$(X,Y)\in\R^2\times\R^2_h$ and $(t,s)\in \R\times\R_{\Delta t}$,
consider the function \beqs\Psi(t,X,s,Y)=\al W^\al(t,X)-\al
W_{h,\Delta
t}^\al(s,Y)+\left(5C_0+\frac{m}{2}\right)\beta_\ep(t-s,X-Y),\eeqs
where, as before, $C_0=\|F(\cdot,\cdot,\cdot,P)\|_\infty$ and
$\beta_\ep=\beta\left(\frac{t}{\ep},\frac{X}{\ep}\right)$ with
$\beta$ a  non-negative smooth function such that
$$\left\{%
\begin{array}{ll}
    \beta(t,X)=1-\left|X\right|^2-\left|t\right|^2, & \hbox{if }
    \left|X\right|^2+\left|t\right|^2\leq \frac{1}{2}, \\
    \beta\leq\frac{1}{2}, & \hbox{if } \frac{1}{2}\leq \left|X\right|^2+\left|t\right|^2\leq 1,\\
    \beta=0, & \hbox{if } \left|X\right|^2+\left|t\right|^2>1. \\
\end{array}%
\right.$$ We have the following lemma:

\begin{lem}\label{crliolem}The function $\Psi$ attains its maximum at a point $(t_0,X_0,s_0, Y_0)$ such that
\begin{itemize}
\item [(i)]
$\Psi(t_0,X_0,s_0, Y_0)\geq 5C_0+\frac{3}{2}m$;
\item [(ii)] $\beta_\ep(t_0-s_0,X_0- Y_0)\geq \frac{3}{5}$.
\end{itemize}
\end{lem}
For the proof, see Lemma 4.1 in \cite{cl}.

Lemma \ref{crliolem} (ii) implies that
$$\beta_\ep(t_0-s_0,X_0- Y_0)=1-\left|\frac{X_0-
Y_0}{\ep}\right|^2-\left|\frac{t_0-s_0}{\ep}\right|^2.$$ Then,
from the inequality $\Psi(s_0,Y_0,s_0,Y_0)\leq
\Psi(t_0,X_0,s_0,Y_0)$ we deduce that \beq\label{t0-s0x0-y0<wal}
\left(5C_0+\frac{m}{2}\right)\left(\left|\frac{X_0-
Y_0}{\ep}\right|^2+\left|\frac{t_0-s_0}{\ep}\right|^2\right)\leq
\al W^\al(t_0,X_0)-\al W^\al(s_0,Y_0)\leq 2C_0.\eeq This implies
that $|t_0-s_0|\rightarrow0$ and $|X_0-Y_0|\rightarrow0$ as
$\ep\rightarrow0$. Moreover, since  $W^\al$ and $W_{h,\Delta
t}^{\al }$ are periodic, we can assume that $(t_0,X_0,s_0,Y_0)$
lies in a compact set of $(\R\times\R^2)^2$. Hence, from
\eqref{t0-s0x0-y0<wal} and the continuity of $W^\al$ we get that
\beq\label{t0-s0x0-y0}\left|\frac{X_0-
Y_0}{\ep}\right|^2+\left|\frac{t_0-s_0}{\ep}\right|^2\rightarrow0\quad
\text{as }\ep\rightarrow0.\eeq

Since $(t_0,X_0)$ is a maximum point of $(t,X)\rightarrow \al
W^\al(t,X)+\left(5C_0+\frac{m}{2}\right)\beta_\ep(t-s_0,X- Y_0)$,
we have \beq\label{Wconteq}\begin{split}&
-\frac{5C_0+\frac{m}{2}}{\al}\partial_t
\beta_\ep(t_0-s_0,X_0-Y_0)+\al
W^\al(t_0,X_0)\\&+F\left(t_0,X_0,-\frac{5C_0+m/2}{\al}D_X\beta_\ep(t_0-s_0,X_0-Y_0)+P\right)\leq
0.\end{split}\eeq

Let $i_0,j_0$ and $n_0$ be such that $X_{i_0,j_0}=Y_0$ and
$s_0=t_{n_0}$. Since $(s_0,Y_0)$ is a minimum point of
$(s,Y)\rightarrow\al W_{h,\Delta t}^{\al
}(s,Y)-\left(5C_0+m/2\right)\beta_\ep(t_0-s,X_0-Y)$, we obtain
$$W_{i_0+1,j_0}^{n_0,\al}-W_{i_0,j_0}^{n_0,\al}\geq
\frac{5C_0+m/2}{\al}[\beta_\ep(t_0-s_0,X_0- Y_0-h
e_1)-\beta_\ep(t_0-s_0,X_0- Y_0)],$$where $e_1=(1,0)^T$.  From the
monotonicity of $g$,
\beq\label{Wineq}\begin{split}&\frac{W_{i_0,j_0}^{n_0,\al}-W_{i_0,j_0}^{n_0-1,\al}}{\Delta
t}+\al
W_{i_0,j_0}^{n_0,\al}+g\left(s_0,Y_0,\frac{5C_0+m/2}{\al}(\Delta_1^+\beta_\ep(t_0-s_0,X_0-\cdot))_{i_0,j_0}+p_x,\right.\\&
\left.(\Delta_1^+W^{n_0,\al})_{i_0-1,j_0}+p_x,(\Delta_2^+W^{n_0,\al})_{i_0,j_0}+p_y,
(\Delta_2^+W^{n_0,\al})_{i_0,j_0-1}+p_y\right)\geq
0.\end{split}\eeq But
$$|(\Delta_1^+\beta_\ep(t_0-s_0,X_0-\cdot))_{i_0,j_0}-e_1\cdot
D_Y\beta_\ep(t_0-s_0,X_0-Y_0)|=\frac{h}{2}|e_1^TD^2_{YY}\beta_\ep(t_0-s_0,X_0-\overline{Y})e_1|,$$
for some $\overline{Y}$ belonging to the segment $(Y_0,Y_0+he_1)$.
Assuming $h$ small enough, so that Lemma \ref{crliolem} (ii)
implies that $|t_0-s_0|^2+|X_0-Y_0|^2+h^2\leq \frac{\ep^2}{2}$, we
obtain that
$D^2_{YY}\beta_\ep(t_0-s_0,X_0-\overline{Y})=\frac{2}{\ep^2}I$,
then
\beq\label{delta-grad}|(\Delta_1^+\beta_\ep(t_0-s_0,X_0-\cdot))_{i_0,j_0}-e_1\cdot
D_Y\beta_\ep(t_0-s_0,X_0-Y_0)|=\frac{h}{\ep^2}.\eeq Now,
\eqref{Wineq}, \eqref{delta-grad} and the monotonicity of $g$
yield
\beqs\begin{split}&\frac{W_{i_0,j_0}^{n_0,\al}-W_{i_0,j_0}^{n_0-1,\al}}{\Delta
t}+\al
W_{i_0,j_0}^{n_0,\al}+g\left(s_0,Y_0,\frac{5C_0+m/2}{\al}e_1\cdot
D_Y\beta_\ep(t_0-s_0,X_0-Y_0)+p_x,\right.\\&
\left.(\Delta_1^+W^{n_0,\al})_{i_0-1,j_0}+p_x,(\Delta_2^+W^{n_0,\al})_{i_0,j_0}+p_y,
(\Delta_2^+W^{n_0,\al})_{i_0,j_0-1}+p_y\right)+\widetilde{C}_1
h\frac{5C_0+m/2}{\ep^2\al}\geq 0.\end{split}\eeqs Repeating
similar estimates for the other arguments in $g$ and for the
derivative with respect to time, we finally find that
\beq\label{finaineqW}\begin{split}&\frac{5C_0+m/2}{\al}\partial_s\beta_\ep(t_0-s_0,X_0-Y_0)+\al
W_{i_0,j_0}^{n_0,\al}+\\& F\left(s_0,Y_0,\frac{5C_0+m/2}{\al}
D_Y\beta_\ep(t_0-s_0,X_0-Y_0)+P\right)+C\frac{h+\Delta
t}{\ep^2\al}\geq 0,\end{split}\eeq where $C$ is independent of
$h,\Delta t,\ep$ and $\al$.

Subtracting \eqref{Wconteq} and \eqref{finaineqW} and using (F2)
we get \beq\label{finalalw-alwh} \al W^\al(t_0,X_0)- \al
W^\al_{h,\Delta t}(s_0,Y_0)\leq C\frac{h+\Delta
t}{\ep^2\al}+\frac{C}{\al}\left|\frac{X_0-Y_0}{\ep}\right|^2+\frac{C}{\al}\left|\frac{t_0-s_0}{\ep}\right|^2,\eeq
where $C$ is independent of $h,\Delta t,\ep$ and $\al$.

Choose $\ep=\ep(\Delta t,h)$ such that $\ep\rightarrow0$ as
$(\Delta t,h)\rightarrow(0,0)$ and $\frac{h+\Delta
t}{\ep^2}\rightarrow0$ as $(\Delta t,h)\rightarrow(0,0)$. From (i)
of Lemma \ref{crliolem} \beqs\begin{split} \sup_{i,j,n}|\al
W^{P,\al}(t_n,x_i,y_j)-\al W_{i,j}^{n,P,\al}|&=m\leq
\sup\Psi-\left(5C_0+\frac{m}{2}\right)\beta_\ep(t_0-s_0,X_0-Y_0)\\&=\al
W^\al(t_0,X_0)- \al W^\al_{h,\Delta t}(s_0,Y_0).\end{split}\eeqs
Then from \eqref{finalalw-alwh} and \eqref{t0-s0x0-y0}, we obtain
\beqs\sup_{i,j,n}|\al W^{P,\al}(t_n,x_i,y_j)-\al
W_{i,j}^{n,P,\al}|\leq \frac{C}{\al}o(1)\quad \text{as }(\Delta
t,h)\rightarrow(0,0).\eeqs From the previous estimate, (ii) of
Lemma \ref{wproplem} and Proposition \ref{wal+fberapprox} we
finally obtain
$$|\overline{F}(P)-\overline{F}_h^{\Delta
t}(P)|\leq \widetilde{K}_1\al+K_1\al+\frac{C}{\al}o(1),$$ and
letting $(h,\Delta t)\rightarrow(0,0)$,  we find that
$$\limsup_{(\Delta t,
h)\rightarrow(0,0)}|\overline{F}(P)-\overline{F}_h^{\Delta
t}(P)|\leq \widetilde{K}_1\al+K_1\al,$$ for any fixed $\al>0$.
This implies that $\lim_{(\Delta t,
h)\rightarrow(0,0)}\overline{F}_h^{\Delta t}(P)=\overline{F}(P)$.
Since $K_1=K_1(P)$ and  $\widetilde{K}_1=\widetilde{K}_1(P)$ are
bounded for $P$
 lying on compact subsets of $\R^{2}$, the convergence is uniform on compact sets.\finedim
\begin{rem}{If $F$ is coercive, then we can get an estimate of the rate of
convergence of $\overline{F}_h^{\Delta t}$ to $\overline{F}$.
Indeed, we have: $$|\overline{F}_h^{\Delta t}-\overline{F}|\leq
(h+\Delta t)^\frac{1}{2},$$see Proposition  A.3 in
\cite{acc}.}\end{rem}

 We conclude this subsection
by recalling the principal properties of $\overline{F}_h^{\Delta
t}$.
\begin{prop}\label{coercivappeffham}Assume (g1)-(g6), (H1)-(H4). Then the approximate effective Hamiltonian $\overline{F}_h^{\Delta
t}$ is Lipschitz continuous with a Lipschitz constant independent
of $h$ and $\Delta t$ and for any $p_x\in\R$
$$\overline{F}_h^{\Delta t}(p_x,0)\geq C_2|p_x|.$$
\end{prop}
\dim For the proof of the Lipschitz continuity of $\overline{F}$,
see the proof of Proposition  A.2 in \cite{acc}.

Let us show the coercivity property. Let $(W_{i,j}^{n,P,\al})$ be
a solution of \eqref{ergodiceq} for $P=(p_x,0)$. Let
$(i_0,j_0,n_0)$ be a maximum point of $(W_{i,j}^{n,P,\al})$, then
\beqs\frac{W_{i_0,j_0}^{n_0,P,\al}-W_{i_0,j_0}^{n_0-1,P,\al}}{\Delta
t}\geq0,\,(\Delta_1^+W^{n_0,P,\al})_{i_0,j_0}\leq0,\,(\Delta_1^+W^{n_0,P,\al})_{i_0-1,j_0}\geq0,\eeqs
\beqs(\Delta_2^+W^{n_0,P,\al})_{i_0,j_0}\leq0,\,
(\Delta_2^+W^{n_0,P,\al})_{i_0,j_0-1}\geq0.\eeqs By the
monotonicity assumption (g1) and \eqref{Fcoercive}, we have
\beqs\overline{F}_h^{\Delta t}(p_x,0)\geq
g(t_{n_0},x_{i_0},y_{i_0},p_x,p_x,0,0)=F(t_{n_0},x_{i_0},y_{i_0},p_x,0)\geq
C_2|p_x|.\eeqs\finedim

\subsection{Long time approximation} A different way to approximate
the effective Hamiltonian is given
by the evolutive Hamilton-Jacobi equation \beq\label{longtime}\left\{%
\begin{array}{ll}
    V_t+F(t,x,y,p_x+D_xV,p_y+D_yV) =0, & (t,x,y)\in(0,+\infty)\times\R^{N+1}, \\
    V(0,x,y)=V_0(x,y), & (x,y)\in\R^{N+1},\\
\end{array}%
\right. \eeq where $V_0$ is bounded and uniformly continuous on
$\R^{N+1}$. Indeed, it is proved in \cite{b} that \eqref{longtime}
admits a unique solution $V$ which is bounded and uniformly
continuous on $[0,T]\times\R^{N+1}$ for any $T>0$, and satisfies
\beqs
\lim_{t\rightarrow+\infty}\frac{V(t,x,y)}{t}=-\overline{F}(P).\eeqs

We approximate \eqref{longtime} by the implicit Eulerian scheme

\beq\label{longtimenum}
\begin{array}{ll}\frac{V_{i,j}^{n+1,P}-V_{i,j}^{n,P}}{\Delta t}+S(t_n,x_i,y_j,h,[V^{n+1,P}]_{i,j})=0\\
V_{i,j}^{0,P}=V_0(x_i,y_j),\end{array}%
\eeq where $S$ is defined as in \eqref{S}. A proof of the
existence of a solution $V=(V_{i,j}^{n,P})$ of \eqref{longtimenum}
is given in \cite{cq} under assumptions (g1)-(g5).

Let $W=(W_{i,j}^{n,P,\al})$ be a solution of \eqref{ergodiceq},
then by comparison, there exist constants $\underline{c}$ and
$\overline{c}$ such that
$$\underline{c}+W_{i,j}^{n,P,\al}-n\overline{F}_h^{\Delta
t}(P)\Delta t\leq V_{i,j}^{n,P}\leq
\overline{c}+W_{i,j}^{n,P,\al}-n\overline{F}_h^{\Delta t}(P)\Delta
t.$$ Since $W$ is bounded, this proves that \beqs
\lim_{n\rightarrow+\infty}\frac{V_{i,j}^{n,P}}{n\Delta
t}=-\overline{F}_h^{\Delta t}(P).\eeqs
\subsection{Approximation of the homogenized problem} We now come
back to the $N$-dimensional homogenized problem \eqref{u0}. From
Theorem \ref{convthmbarles} we know that if $\overline{H}$ is the
effective Hamiltonian in \eqref{u0}, then
$\overline{H}(p)=\overline{F}(p,-1)$ for any $p\in\R^N$. Hence,
from Proposition \ref{convappham}, the discrete Hamiltonian
$$\overline{H}_h^{\Delta
t}(p):=\overline{F}_h^{\Delta t}(p,-1),$$ is an approximation of
$\overline{H}(p)$ for any $p\in\R^N$.

As in \cite{acc}, we approximate \eqref{u0} by the problem
\beq\label{u0app} \left\{%
\begin{array}{ll}
    \partial_t u_{\Delta t,h}+\overline{H}_{h}^{\Delta t}(Du_{\Delta t,h})=0, & (t,x)\in(0,+\infty)\times\R^N, \\
    u_{\Delta t,h}(0,x)=u_0(x), & x\in\R^N,\\
\end{array}%
\right. \eeq where $h$ and $\Delta t$ are fixed, and $u_0$ is the
same initial datum as in \eqref{u0}.

By Proposition \ref{coercivappeffham} $\overline{H}_{h}^{\Delta
t}$ is Lipschitz continuous and coercive, so \eqref{u0app} has a
unique viscosity solution $u_{\Delta t,h}$ which is an
approximation of the solution $u^0$ of \eqref{u0}:
\begin{prop}\label{u0appudh} Let $u^0$ and $u_{\Delta t,h}$ be respectively the viscosity solutions of \eqref{u0} and \eqref{u0app}. Then
for any $T>0$
\beq\label{u0hconvu0}\sup_{[0,T]\times\R^N}|u_{\Delta
t,h}-u^0|\rightarrow0\quad \text{as }(\Delta t,
h)\rightarrow(0,0).\eeq
\end{prop}
\dim If $L_0$ is the Lipschitz constant of the initial datum
$u_0$, then, by Proposition \ref{regulU0}, the functions
$u^0(t,\cdot)$ and $u_{\Delta t,h}(t,\cdot)$ are Lipschitz
continuous  with same Lipschitz constant $L_0$. By Proposition
\ref{convappham} the approximate Hamiltonian
$\overline{H}_{h}^{\Delta t}$ converges to $\overline{H}$
uniformly for $|p|\leq L$. Hence \eqref{u0hconvu0} follows by the
following proposition, which is a standard estimate in the regular
perturbation theory of Hamilton-Jacobi equations
 (see Theorem VI.22.1 in \cite{bc})
\begin{prop}If there exists $\eta>0$ such that if $H_i$, $i=1,2$, satisfy (H1)-(H3) with $$\|H_1-H_2\|_\infty\leq\eta,$$ and if $u_i$, $i=1,2$,
are viscosity solutions of \beqs\left\{
       \begin{array}{ll}
         u_t+H_i(Du)=0, & (t,x)\in(0,T)\times\R^N \\
         u(0,x)=u_0(x), & x\in\R^N,
\end{array}
     \right.
\eeqs where $u_0$ is bounded and uniformly continuous on $\R^N$,
then, for some constant $C$,
$$\|u_1-u_2\|_\infty\leq C\eta.$$

\end{prop}
\finedim
\begin{rem}{\em In order to compute numerically the approximation of $u^0$, we need  further discretizations. Indeed, we have approximated
$\overline{H}(p)$ by $\overline{H}_{h}^{\Delta t}(p)$ for any
fixed $p\in\R^N$. Since it is not possible to compute
$\overline{H}_h^{\Delta t}(p)$ for any $p$, one possibility is to
introduce a triangulation of a bounded region of $\R^N$ and compute
$\overline{H}_h^{\Delta t}(p_i)$, where $p_i$ are the vertices of
the simplices and to approximate all the other values
$\overline{H}_h^{\Delta t}(p)$ by  $\overline{H}_{h,k}^{\Delta
t}(p)$, where $\overline{H}_{h,k}^{\Delta t}$ is the linear
interpolation of $\overline{H}_h^{\Delta t}$ and  we denote by $k$
the maximal diameter of the simplices. The solution $u_{\Delta
t,h}^k$ of
\beq\label{interphompb} \left\{%
\begin{array}{ll}
    \partial_t u_{\Delta t,h}^k+\overline{H}_{h,k}^{\Delta t}(Du_{\Delta t,h}^k)=0, & (t,x)\in(0,+\infty)\times\R^N, \\
    u_{\Delta t,h}^k(0,x)=u_0(x), & x\in\R^N,\\
\end{array}%
\right. \eeq is an approximation of $u_{\Delta t,h}$ as
$k\rightarrow0$ and hence, by Proposition \ref{u0appudh}, of $u^0$
as $(\Delta t, h, k)\rightarrow(0,0,0)$. Finally, discretizing
\eqref{interphompb} by means a monotone, consistent and stable
approximation scheme, we can compute numerically an approximation
of the solution $u^0$ of \ref{u0}. See \cite{acc} for details.
}\end{rem}

\newcommand{\dt}{\Delta t}
\newcommand{\shape}{\mathcal{Q}}
\newcommand{\HH}{{\mathbb H}}
\newcommand{\V}{{\mathbb V}}

\newcommand{\RN}{{\mathbb R^N}}

\newcommand{\zz}{{\mathbb z}}

\newcommand{\bO}{{\mathbb O}}
\newcommand{\bM}{{\mathbb M}}
\newcommand{\bL}{{\mathbb L}}

\newcommand{\cF}{{\mathcal F}}
\newcommand{\cG}{{\mathcal G}}
\newcommand{\cE}{{\mathcal E}}
\newcommand{\cI}{{\mathcal I}}
\newcommand{\cJ}{{\mathcal J}}
\newcommand{\cA}{{\mathcal A}}
\newcommand{\cN}{{\mathcal N}}
\newcommand{\cS}{{\mathcal S}}
\newcommand{\cL}{{\mathcal L}}
\newcommand{\cD}{{\mathcal D}}
\newcommand{\cC}{{\mathcal C}}
\newcommand{\cH}{{\mathcal H}}
\newcommand{\cM}{{\mathcal M}}
\newcommand{\cP}{{\mathcal P}}
\newcommand{\cU}{{\mathcal U}}
\newcommand{\cV}{{\mathcal V}}
\newcommand{\cT}{{\mathcal T}}
\newcommand{\cO}{{\mathcal O}}
\newcommand{\cR}{{\mathcal R}}

\newcommand{\Exp}{{^}}
\newcommand{\deter}{\hbox{det}}
\newcommand{\Ker}{\ker}
\newcommand{\conv}{\hbox{conv}}

\newcommand{\tv}{{\tilde v}}
\newcommand{\ds}{\displaystyle}
\newcommand{\supp}{\mathrm{supp}\;}

\section{Numerical Tests}
\label{sec:numerical-tests}
The present paragraph is devoted to the description of numerical approximations of the effective Hamiltonian.
 \subsection{Results}
\subsubsection{First case}
\begin{figure}[htbp]
  \centering
 \[ \begin{array}[c]{l}
    \includegraphics[width=12cm]{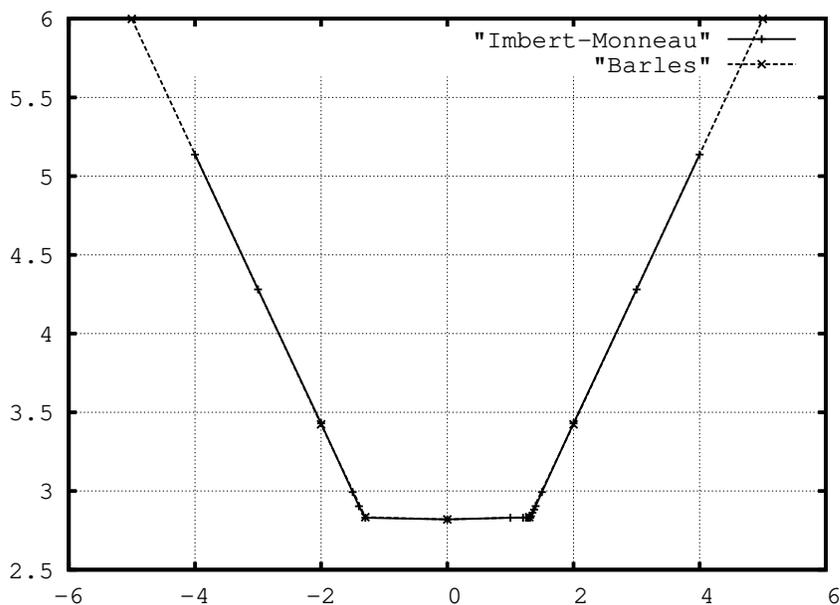}
  \end{array}\]
  \caption{First case: the effective Hamiltonian as a function of $p$ obtained with both Barles and Imbert-Monneau cell problems.}
  \label{fig:1}
\end{figure}

\begin{figure}[htbp]
  \centering
 \[ \begin{array}[c]{l}
 \includegraphics[width=8cm]{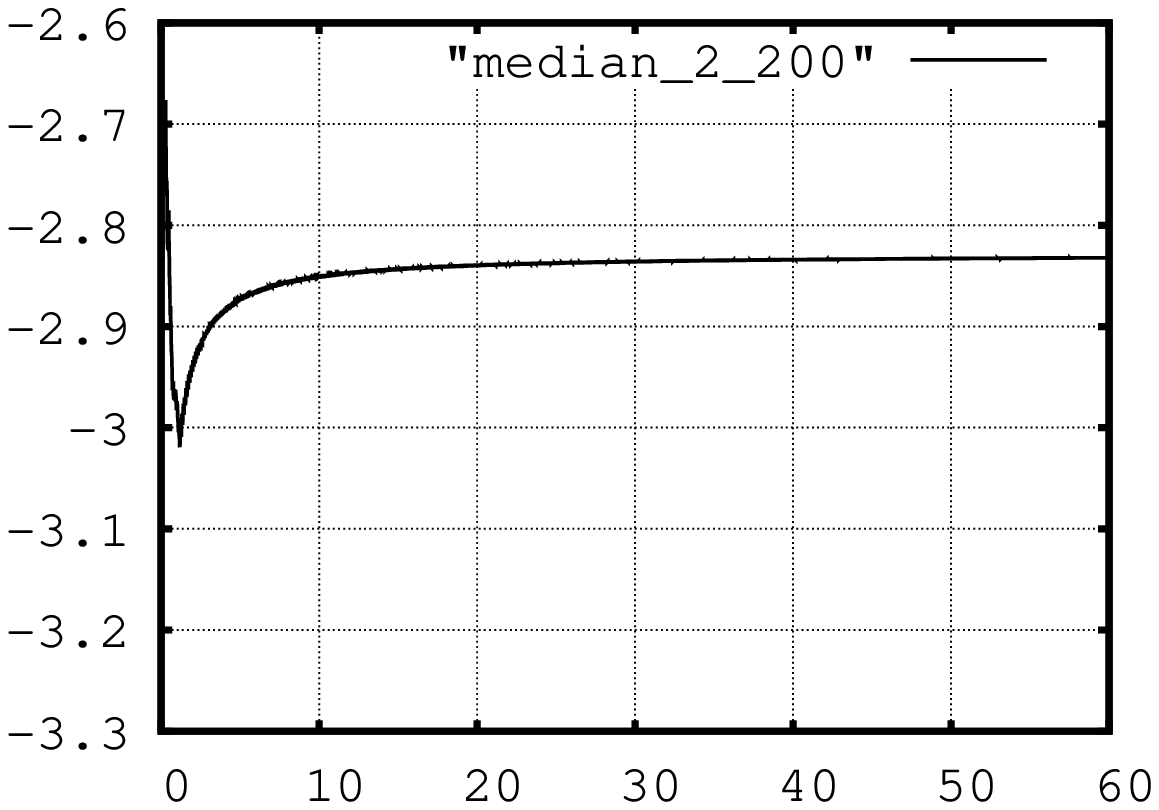}
    \includegraphics[width=8cm]{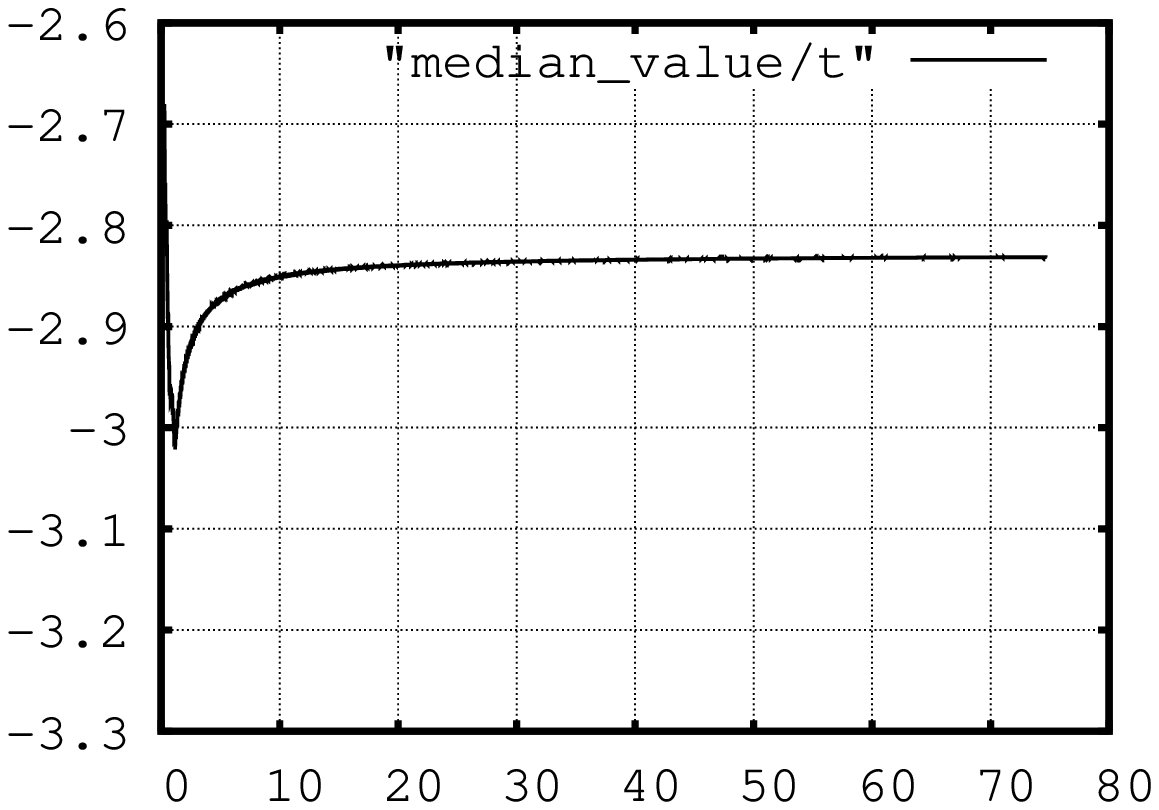}
  \end{array}\]
  \caption{First case, $p=1.3$. Left: the median value  of $w(\tau,\cdot)/\tau$ on a period as a function of $\tau$. Right: the median value  of $v(\tau,\cdot)/\tau$ on a period as a function of $\tau$}
  \label{fig:2}
\end{figure}

\begin{figure}[htbp]
  \centering
 \[ \begin{array}[c]{l}
    \includegraphics[width=8cm]{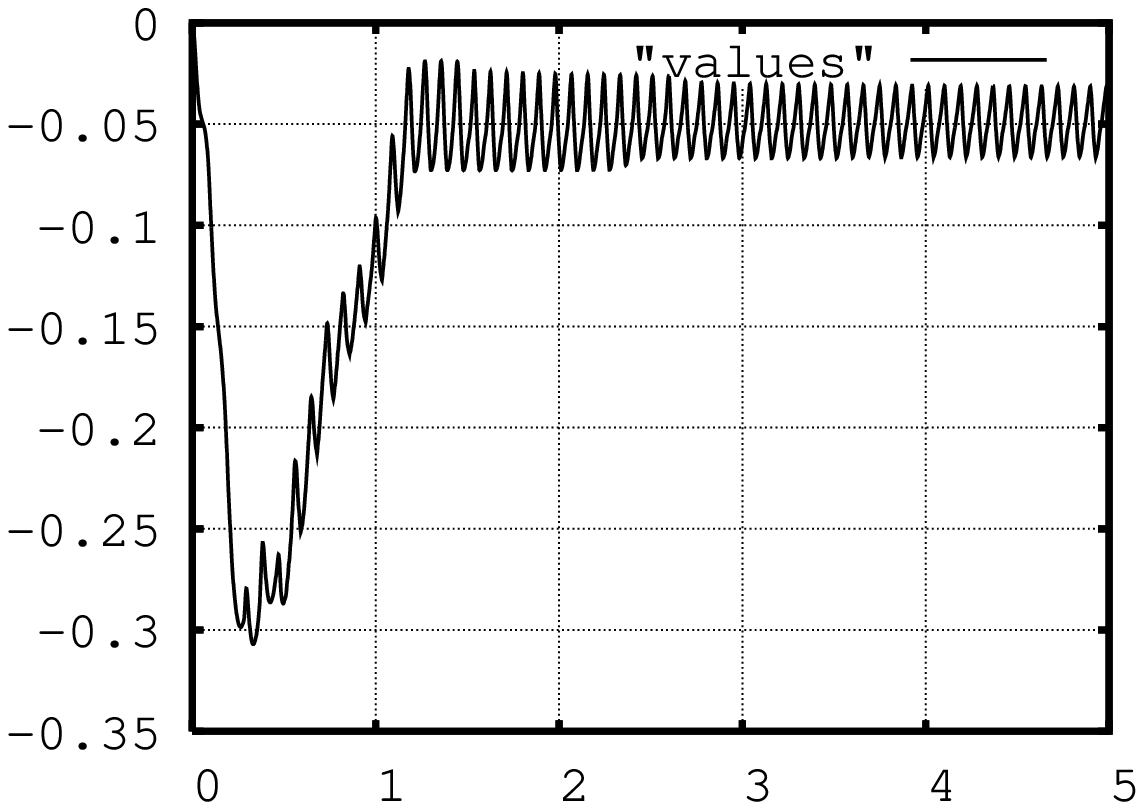}
    \includegraphics[width=8cm]{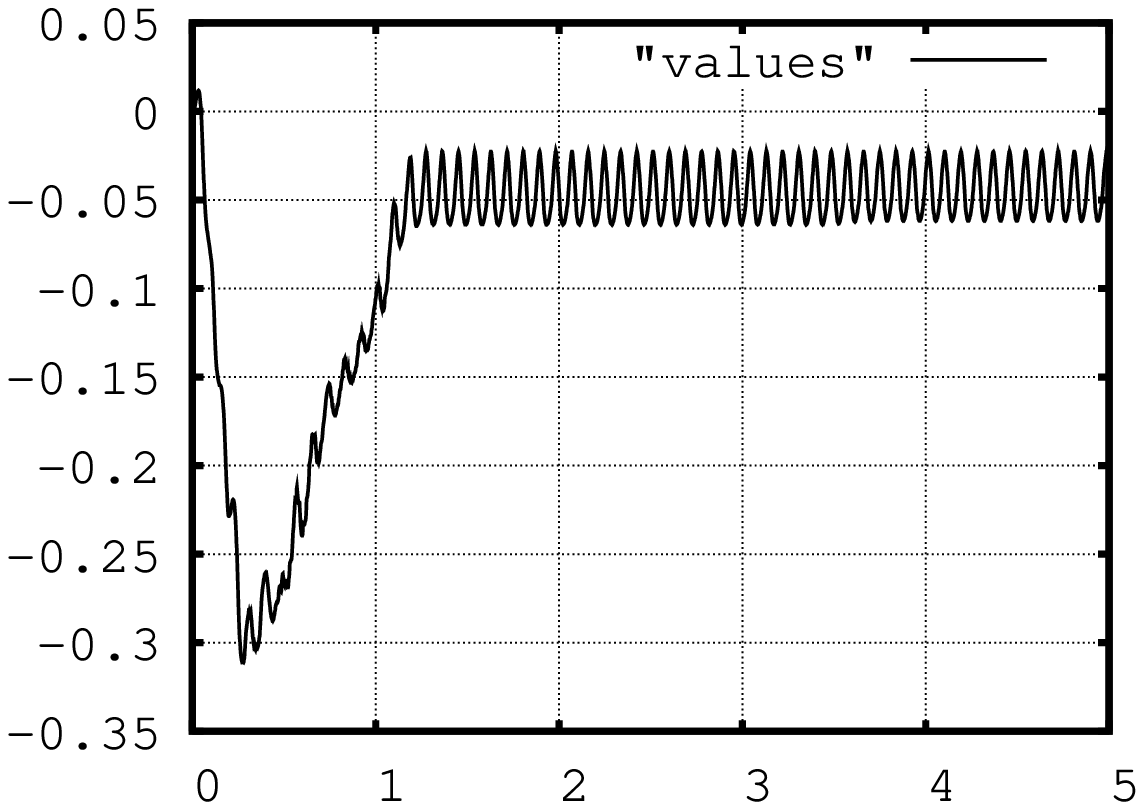}
  \end{array}\]
  \caption{First case, $p=1.3$: $w(\tau,0,0)-\langle w(\tau)\rangle$ (left) and  $v(\tau,0)-\langle v(\tau)\rangle$ (right) as a function of $\tau$}
  \label{fig:3}
\end{figure}

\begin{figure}[htbp]
  \centering
 \includegraphics[width=14cm]{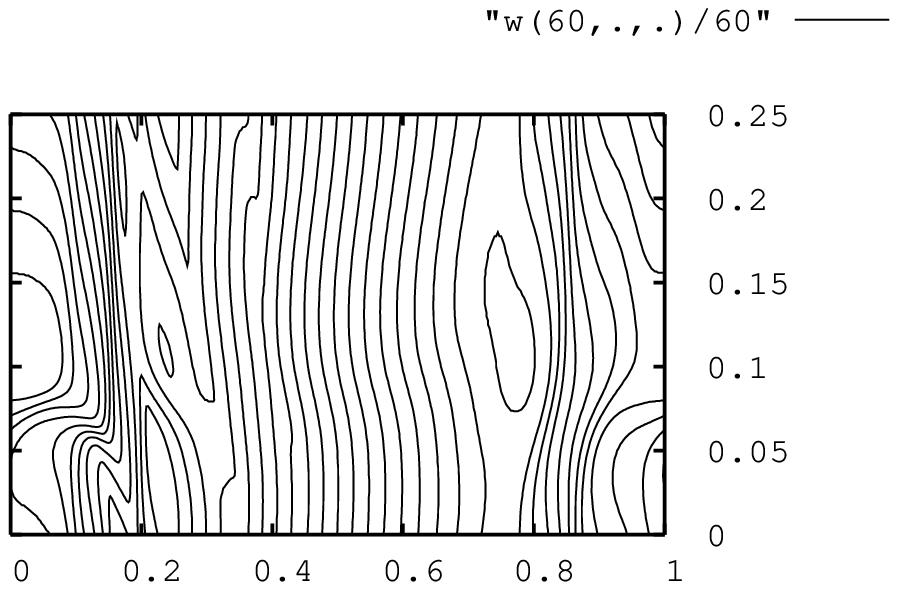}
  \includegraphics[width=12cm]{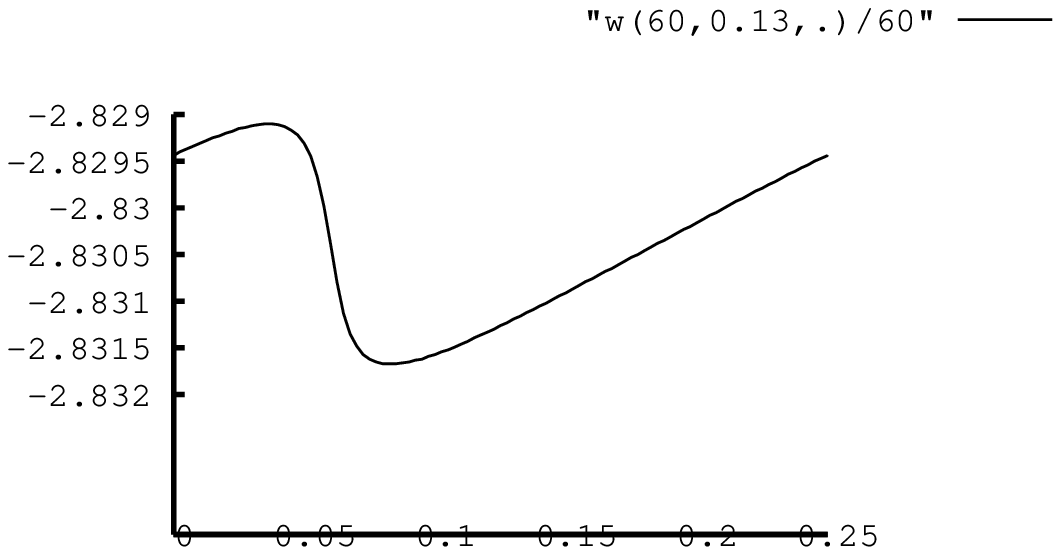}
  \caption{First case, Barles cell problem, $p=1.3$. Top: contour lines of $w(\tau,\cdot)/\tau$ on a period as a function of $(x,y)$. Bottom: the cross-section  $x=0.13$.}
  \label{fig:4}
\end{figure}

\begin{figure}[htbp]
  \centering
 \[ \begin{array}[c]{ll}
 \!\!\!\!\!\!\!\!\!\!\!\!\!\!\!\!\!\!\!\!\!
    \includegraphics[width=10cm]{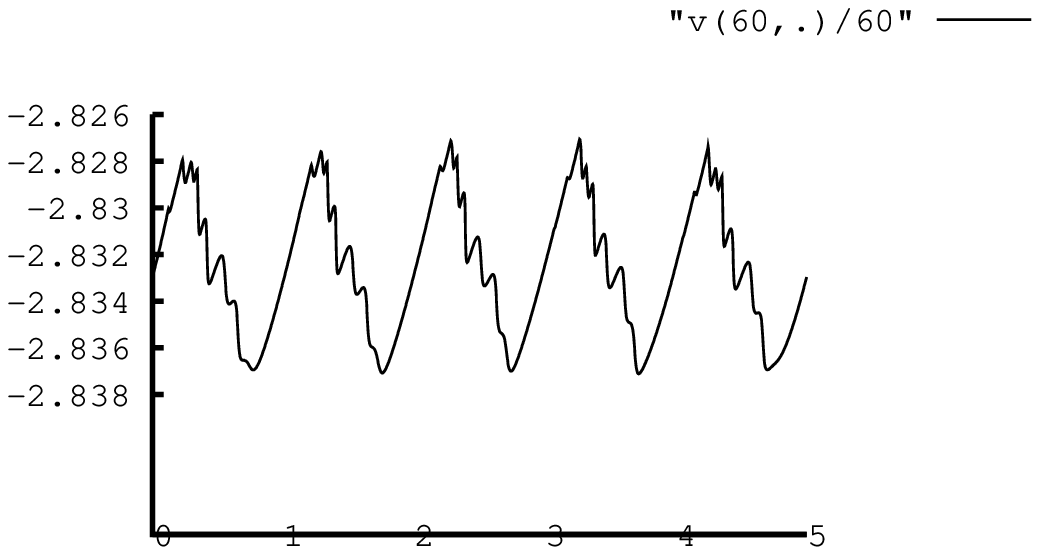}  \!\!\!\!\!\!\!\!\!\!\!\!\!\!\!\!\!\!\!\!\!\!\!\!\!\!\!\!\!\!\!\!\!\!\!\!\!\!\!\!   &\includegraphics[width=10cm]{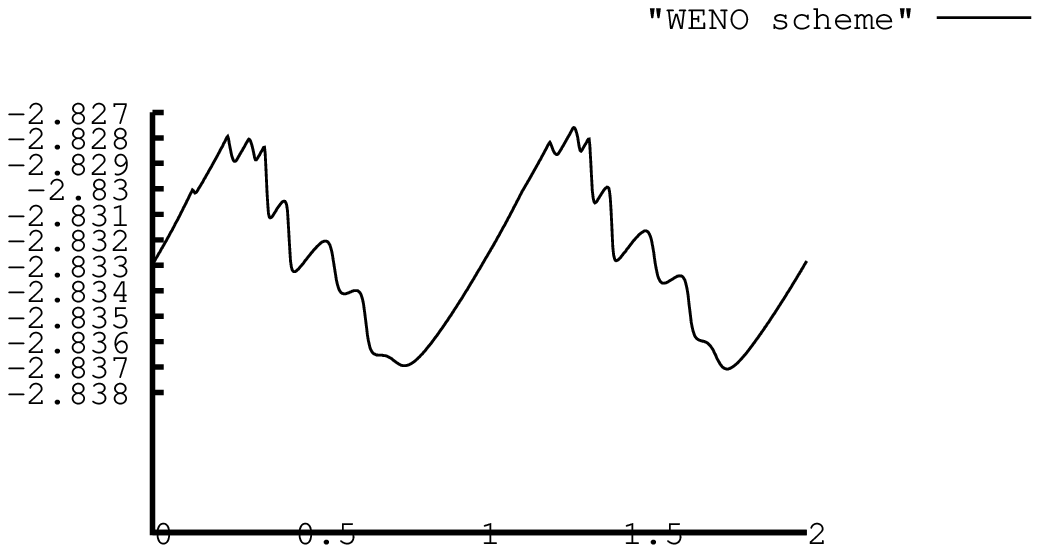}
\\ &\includegraphics[width=10cm]{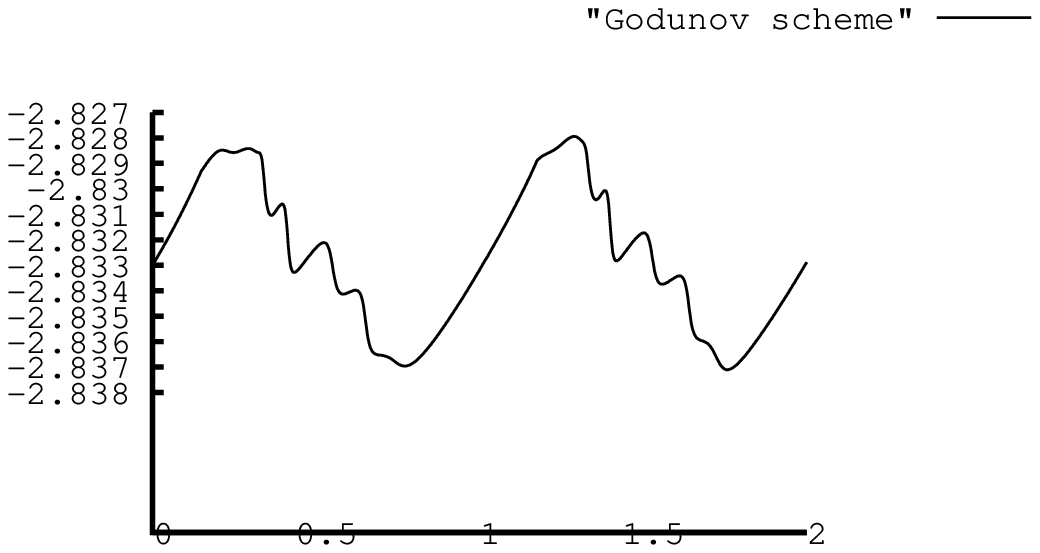}
  \end{array}\]
  \caption{First case, Imbert-Monneau cell problem, $p=1.3$: Top: Third order Runge Kutta/WENO scheme: $v(\tau,x)/\tau$ as a function of $x$ for $\tau=60$; the right part is a zoom. Bottom: same computation with Euler/Godunov scheme with the same grid parameters: some oscillations are smeared out, but the average value of the solution is well computed.}
  \label{fig:5}
\end{figure}

 We discuss a one dimensional case where the Hamiltonian is
 \begin{displaymath}
   H(x,u,p)= 2\cos(2 \pi x) +\sin(8\pi u) + (1-\cos(6\pi x) /2) |p|.
 \end{displaymath}
We have used two approaches for computing the effective Hamiltonian.
\begin{enumerate}
\item Barles cell problem: the first approach consists of increasing the dimension and considering the
long time behavior of the continuous viscosity solution $w$ of
\begin{equation}\label{eq:num_2}
   \begin{array}[c]{rcl}
\ds   w_t +F(x,y, p+D_x w,-1+D_yw  )=0, \quad &&\ds(t,x,y)\in (0,\infty)\times \R\times \R,\\
   \ds w(0,x,y)=0, &&(x,y)\in\R\times \R,
   \end{array}
\end{equation}
where $F$ is given by (\ref{F}).
In the present case, from the periodicity of $H$ with respect to $x$ and $u$, $w$ is $1$-periodic with respect to $x$ and
$1/4$-periodic with respect to $y$.
 We know that when $t\to \infty$,  $w(t,\cdot,\cdot)/t$ tends to a real number $\lambda$ and that
$\overline H(p)=-\lambda$. \\
For approximating (\ref{eq:num_2}) on a uniform grid, we have used an
explicit Euler time marching method with a Godunov monotone scheme
(see \cite{MR606500,MR1700751}).
 A semi-implicit time marching scheme which allows for large time steps may be used as well, see \cite{acc},
 but very large time steps cannot be taken because of the periodic in time asymptotic behaviour of $w$.
\\
Alternatively, we have also used the higher order method described in \cite{MR1762034}, see also \cite{MR1391627}.
It is a third order TVD explicit Runge-Kutta time marching method with a weighted ENO scheme
 in the spatial variables.
This weighted ENO scheme is constructed upon and has the same stencil nodes
 as the third order ENO scheme but  can be as as high as fifth order accurate in the smooth part of the solution.
\item Imbert-Monneau cell problem:
when $p$ is a rational number ($p= \frac {n}{q}$),
instead of considering a problem posed in two space dimensions,  one possible way of approximating the
 effective Hamiltonian $\overline H(p)$ is to consider the cell problem
 \begin{equation}\label{eq:num_1}
   \begin{array}[c]{rcl}
\ds   v_t +H(x,v+p\cdot x, p+Dv)=0, \quad &&\ds(t,x)\in (0,\infty)\times \R,\\
   \ds v(0,x)=0 &&x\in\R.
   \end{array}
 \end{equation}
This problem has a unique continuous solution which is periodic of
period $q$ with respect to $x$ (in fact, the smallest period of
$v$ may be a divisor of $q$). From \cite{im} (Theorem 1), we know
that there exists a unique real number $\lambda$ such that $\frac
{v(\tau,x)}{\tau}$ converges to $\lambda$ as $\tau\to \infty$
uniformly in $x$, and that $\overline H(p)=-\lambda$. Moreover,
when $t$ is large, the function $v(t,x)-\lambda t$ becomes close
to a periodic function of time. In what follows, (\ref{eq:num_1})
will be referred to as Imbert-Monneau cell problem. Note that the
size of the period varies with $p$ and may be arbitrary large.
This is clearly a drawback of this approach which is yet the
fastest one for one dimensional problems and
 moderate values of $q$.\\
For approximating (\ref{eq:num_1}) on a uniform grid, we have
used either the abovementioned explicit Euler time marching method with a Godunov monotone scheme
or the third order TVD explicit Runge-Kutta time marching method with a weighted ENO scheme
 in the spatial variable.
\end{enumerate}
In Figure~\ref{fig:1}, we plot the graph of the effective Hamiltonian computed with the high order methods and both Imbert-Monneau and Barles cell problems.
For Barles cell problems, the grid of the square $[0,1]\times[0,1/4]$ has $400\times 100$ nodes and the time step is $1/1000$.
For Imbert-Monneau cell problems, the grids in the $x$ variable are uniform with a step of $1/400$ and the time step is $1/1000$.
The two graphs are undistinguishable.
It can be seen that the effective Hamiltonian is symmetric with respect to $p$ and constant for small values of $p$, i.e. $|p|\lesssim 1.3$.
The points where we have computed the effective Hamiltonian are concentrated near $1.3$ where the slope of the graph changes.
Our computations clearly indicate that the effective Hamiltonian is piecewise linear.
\\
In order to show the convergence of $\frac {v(\tau,x)}{\tau}$ and  $\frac {w(\tau,x)}{\tau}$, we take $p=1.3$
so the space period of the Imbert-Monneau cell problem is $5$.
In Figure ~\ref{fig:2}, we plot  $\frac {\langle w(\tau)\rangle} \tau$ (left) and   $\frac {\langle v(\tau)\rangle} \tau$ (right)
as a function of $\tau$,
where  $\langle v(\tau)\rangle$ is the median value of $v(\tau,\cdot)$ on a spatial period. Both functions converge to constants when $\tau\to \infty$ and the limit are close to each other (the error between the two scaled median values is smaller than  $10^{-3}$ at $\tau \sim 60$ and we did not consider much longer times).
In Figure~\ref{fig:3}, we plot the graphs of the functions $w(\tau,0,0)-\langle w(\tau)\rangle$ (left) and
$v(\tau,0)-\langle v(\tau)\rangle$ (right).
 We see that these functions become close to  time-periodic.
In Figure~\ref{fig:4} (top), we plot the contour lines of the function $w(\tau,x,y)/\tau$
as a function of $(x,y)$ for $\tau=60$. In the bottom part of the figure we plot the graph of $y \to  w(\tau,0.13,y)/\tau$
for the same value of $\tau$. We see that $w$ has internal layers.
In Figure~\ref{fig:5}, we plot the graph $x\to v(\tau,x)/\tau$  for $\tau=60$.
We first see that the function takes all its values in a small interval and  has
 very rapid variations with respect to $x$ (is nearly discontinuous). This does not contradict the theory, because there are
no uniform estimates on the modulus of continuity of $v(\tau,\cdot)/\tau$.

\subsubsection{Second case}
\begin{figure}[htbp]
  \centering
 \[ \begin{array}[c]{l}
    \includegraphics[width=12cm]{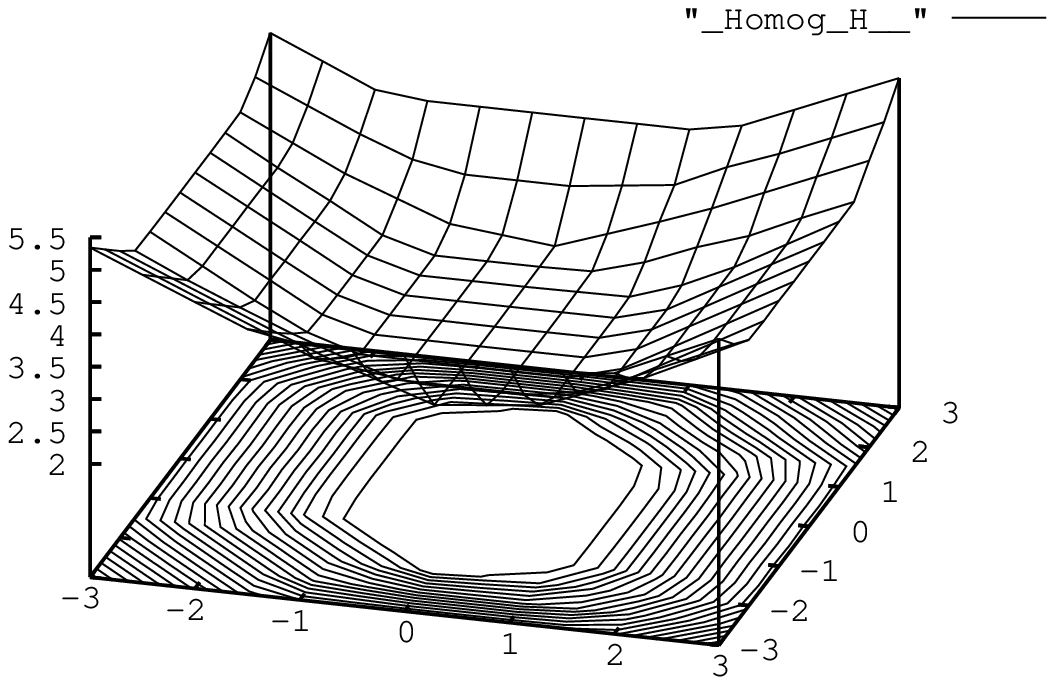}
  \end{array}\]
  \caption{Second case, the effective Hamiltonian computed by solving Imbert-Monneau cell problems.}
  \label{fig:6}
\end{figure}

\begin{figure}[htbp]
  \centering
 \[ \begin{array}[c]{l}
 \includegraphics[width=8cm]{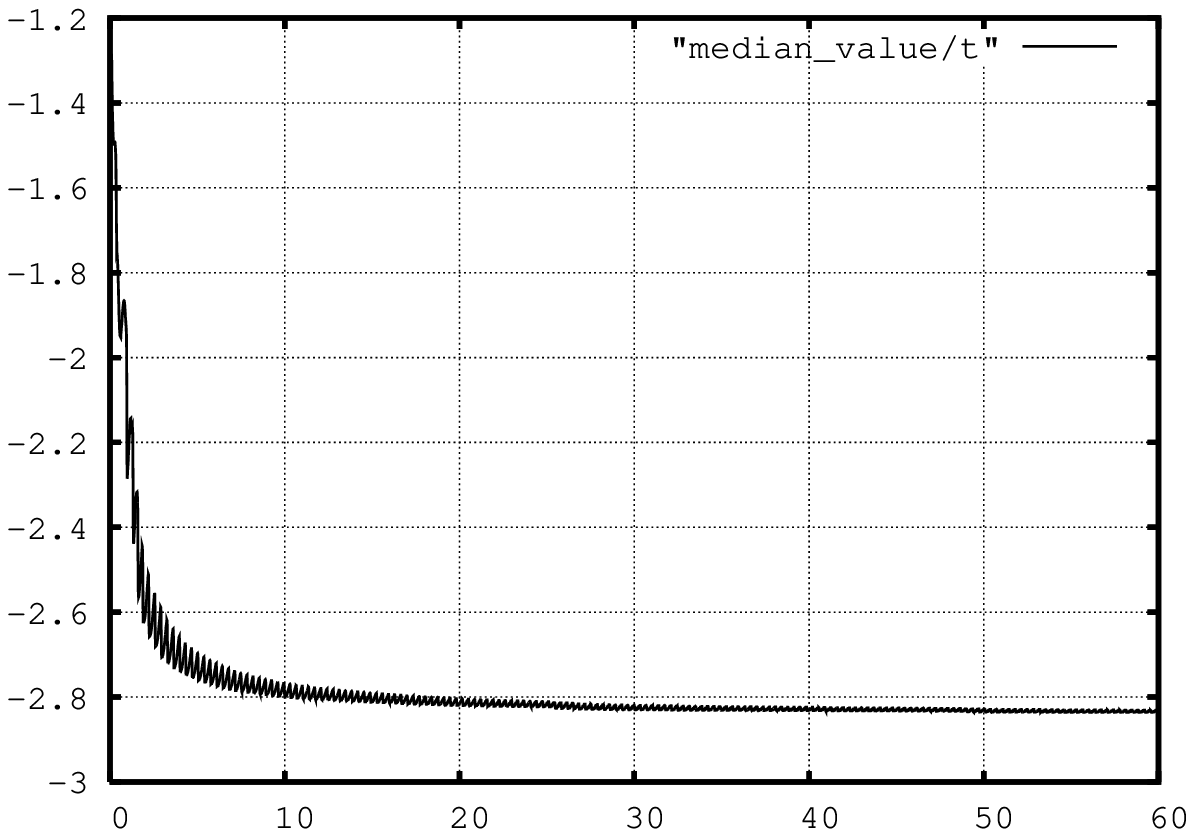}
    \end{array}\]
  \caption{Second case, $p=(1,1)$. The median value  of $v(\tau,\cdot)/\tau$ on a period as a function of $\tau$. }
  \label{fig:7}
\end{figure}

\begin{figure}[htbp]
  \centering
 \[ \begin{array}[c]{l}
    \includegraphics[width=14cm]{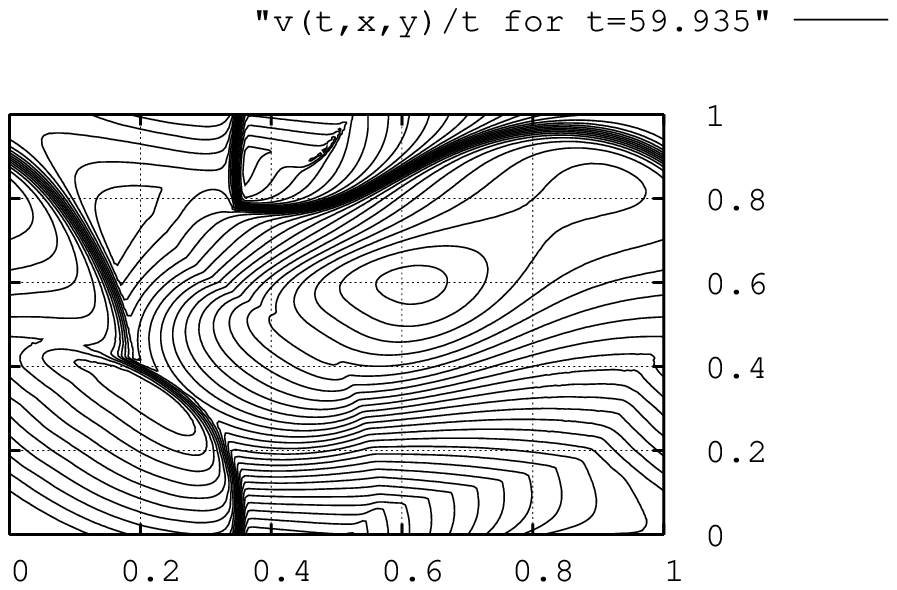}
  \end{array}\]
  \caption{Second case, the contours of the solution of Imbert-Monneau cell problem for $p=(1,1)$ at time $\tau=59.935$.}
  \label{fig:8}
\end{figure}

We consider a two dimensional problem, where the Hamiltonian is
 \begin{displaymath}
   H(x,u,p)= \cos(2\pi x_1)+\cos(2\pi x_2)+\cos(2\pi(x_1-x_2)) +\sin(2\pi u) +
\left(1-\frac {\cos(2\pi x_1)}2-\frac {\sin(2\pi x_2)}4\right) |p|.
 \end{displaymath}
For this case, only the Imbert-Monneau cell problems have been approximated on uniform grids with step $1/200$. The time step is $0.005$.
 In Figure~\ref{fig:6}, we plot the contours and the graph of the effective Hamiltonian computed with the high order method.
We can see that the effective Hamiltonian is symmetric with respect to $p=(0,0)$,  constant for small vectors $p$.
In Figure ~\ref{fig:7}, we plot   $\frac {\langle v(\tau)\rangle} \tau$
as a function of $\tau$. We see that this function converges when $\tau\to \infty$.
In Figure ~\ref{fig:8}, we plot the contours of $v(\tau,\cdot)/\tau$ for $\tau=59.935$ and $p=(1,1)$.
We see that for large values of $\tau$, $v$ is close to discontinuous.

% \bibliographystyle{plain}
% \bibliography{numerical_test}

%\end{document}

\appendix
\section{}
\dims {\bf of Lemma \ref{varphilem}.} To show that the sequence is
convergent it suffices to show that for any $s\in\R$
$\varphi_\ep^{n,\delta}(s)$ is a Cauchy sequence. Fix $s\in\R$ and
let $i_0\in\Z$ be the closest integer to $s$, i.e., $s=
i_0\ep+\gamma\ep$, with
$\gamma\in\left(-\frac{1}{2},\frac{1}{2}\right]$. Let $k>m>|i_0|$,
then, by assumptions \eqref{varphiprop} we have
\beqs\begin{split}\varphi_\ep^{k,\delta}(s)-\varphi_\ep^{m,\delta}(s)&=\sum_{i=-k}^{-m-1}\ep\phi\left(\frac{s-\ep
i}{\delta}\right)+\sum_{i=m+1}^{k}\ep\phi\left(\frac{s-\ep
i}{\delta}\right)-\ep(k-m)\\&=\sum_{i=-k}^{-m-1}\ep\left[\phi\left(\frac{s-\ep
i}{\delta}\right)-1\right]+\sum_{i=m+1}^{k}\ep\phi\left(\frac{s-\ep
i}{\delta}\right)\\& \leq \ep
K_2\delta^2\sum_{i=-k}^{-m-1}\frac{1}{(s-\ep i)^2}+\ep
K_2\delta^2\sum_{i=m+1}^{k}\frac{1}{(s-\ep i)^2}\\&
=K_2\frac{\delta^2}{\ep}\sum_{i=-k}^{-m-1}\frac{1}{(i_0-i+\gamma)^2}+K_2\frac{\delta^2}{\ep}\sum_{i=m+1}^{k}\frac{1}{(i_0-i+\gamma)^2}.
\end{split}\eeqs

Similarly, it can be showed that
$$\varphi_\ep^{k,\delta}(s)-\varphi_\ep^{m,\delta}(s)\geq -K_2\frac{\delta^2}{\ep}\sum_{i=-k}^{-m-1}\frac{1}{(i_0-i+\gamma)^2}
-K_2\frac{\delta^2}{\ep}\sum_{i=m+1}^{k}\frac{1}{(i_0-i+\gamma)^2}.$$
Hence
$|\varphi_\ep^{k,\delta}(s)-\varphi_\ep^{m,\delta}(s)|\rightarrow0$
as $m,k\rightarrow+\infty$. Similar arguments show that the
sequence $(\varphi_\ep^{\delta,n})'$ converge uniformly on compact
sets of $\R$. This implies that $\varphi_\ep^\delta$ is of class
$C^1$ with
$(\varphi_\ep^\delta)'(s)=\lim_{n\rightarrow+\infty}(\varphi_\ep^{\delta,n})'(s)$.

Now, let us show \eqref{varphilim}. Let $s=i_0\ep+\gamma\ep$ for
some $i_0\in\Z$ and $\gamma\in[0,1)$. Then

\beqs\begin{split}\varphi_\ep^{n,\delta}(s)-i_0\ep &
=\ep\left[\phi\left(\frac{\gamma\ep}{\delta}\right)-1\right]+\sum_{i=-n}^{i_0-1}\ep\left[\phi\left(\frac{i_0\ep+\gamma\ep-\ep
i}{\delta}\right)-1\right]+\sum_{i=i_0+1}^{n}\ep\phi\left(\frac{i_0\ep+\gamma\ep-\ep
i}{\delta}\right)\\&\leq
\ep\left[\phi\left(\frac{\gamma\ep}{\delta}\right)-1\right]\ep+\frac{\delta^2}{\ep}K_2\sum_{i=-n}^{i_0-1}\frac{1}{(i_0-i+\gamma)^2}
+\frac{\delta^2}{\ep}K_2\sum_{i=i_0+1}^{n}\frac{1}{(i-i_0-\gamma)^2}\\&
=\ep\left[\phi\left(\frac{\gamma\ep}{\delta}\right)-1\right]+\frac{\delta^2}{\ep}K_2\sum_{i=1}^{n+i_0}\frac{1}{(i+\gamma)^2}
+\frac{\delta^2}{\ep}K_2\sum_{i=1}^{n-i_0}\frac{1}{(i-\gamma)^2}.
\end{split}\eeqs
Similarly
$$\varphi_\ep^{n,\delta}(s)-i_0\ep\geq \ep\left[\phi\left(\frac{\gamma\ep}{\delta}\right)-1\right]-
\frac{\delta^2}{\ep}K_2\sum_{i=1}^{n+i_0}\frac{1}{(i+\gamma)^2}
-\frac{\delta^2}{\ep}K_2\sum_{i=1}^{n-i_0}\frac{1}{(i-\gamma)^2}.$$
Letting $n\rightarrow+\infty$, we get
$$\left|\varphi_\ep^{\delta}(s)-i_0\ep-\ep\left[\phi\left(\frac{\gamma\ep}{\delta}\right)-1\right]\right|\leq
\frac{\delta^2}{\ep}K_2\sum_{i=1}^{+\infty}\frac{1}{(i+\gamma)^2}
+\frac{\delta^2}{\ep}K_2\sum_{i=1}^{+\infty}\frac{1}{(i-\gamma)^2}.$$
If $\gamma>0$ then
$\phi\left(\frac{\gamma\ep}{\delta}\right)-1\rightarrow0$ as
$\delta\rightarrow0^+$ and $\varphi_\ep^{\delta}(s)\rightarrow
i_0\ep$ if $\delta\rightarrow0^+$. If $\gamma=0$, then
$\varphi_\ep^{\delta}(s)\rightarrow (i_0-1)\ep+\phi(0)\ep$ if
$\delta\rightarrow0^+$ and \eqref{varphilim} is proved.\finedim
\begin{ack}{\em The second author was supported by DEASE:
MEST-CT-2005-021122.}\end{ack}

\end{document}